\documentclass[10pt,twoside]{article}
\usepackage{amssymb}
\usepackage[mathscr]{eucal}
\usepackage{eufrak}
\usepackage{amsmath}
\usepackage{mathrsfs}
\font\bg=cmbx10 scaled 1200

\def\intll#1#2{\int\limits_{#1}^{#2}}

\def\rf#1{{\rm(\ref{#1})}}
\def\chiu{\hfill$\displaystyle\vspace{4pt}
\underset\Box\null$}
\def\Pr{{\bf Proof. }}
\def\O{\Omega}

\def\R{\Bbb R}
\def\N{\Bbb N}
\def\o{\"{o}}
\def\à{\`{a}}
\def\è{\`{e}}
\def\ì{\`{i}}
\def\ù{\`{u}}
\def\ò{\`{o}}
\def\é{\'{e}}
\def\vf{\varphi}
\def\dy{\displaystyle}
\def\ve{\varepsilon}
\def\l{\lambda}
\def\pa{\partial}
\def\be{\begin{equation}}
\def\ba{\begin{array}}
\def\ea{\end{array}}
\def\ee{\end{equation}}
\def\vs1{\vspace{1ex}}
\def\vp{\varphi}
\def\ov{\overline}
\def\po{\partial\Omega}
\font\sc=cmcsc10
\title{\bg Higher regularity of
solutions\\ to  the singular
\mbox{\large$p\,$}-Laplacean
parabolic system}
\author{\sc F. Crispo and P. Maremonti
\thanks{
Dipartimento di Matematica, Seconda
Universit\`{a} degli Studi di
 Napoli, via Vivaldi 43, 81100 Caserta,
 Italy.
francesca.crispo@unina2.it; paolo.maremonti@unina2.it}}
\date{\small{2012/08/30}}
\begin{document}
\maketitle
\noindent{\bf Abstract}  - {\small We study
 existence and regularity properties of solutions to the
singular $p$-Laplacean parabolic system in a bounded
 domain $\O$. The main purpose is
to prove global
 $L^r(\ve,T;L^q(\O))$, $\ve\geq0$, integrability
 properties of the second spatial
derivatives and of the time
derivative of the solutions. Hence, for suitable
$p$ and exponents $r,\,q$,
by Sobolev embedding theorems, we
deduce global regularity
of $u$ and $\nabla u$ in
H\"older spaces. Finally we prove a global pointwise bound for the
solution under the assumption $p>\frac{2n}{n+2}$. }\vskip 0.2cm
 \par\noindent{Keywords: parabolic
system, singular $p$-Laplacean, higher integrability,
 global regularity.}
 \par\noindent
{\par\noindent M.R.: }
 \vskip -0.7true cm\noindent
\newcommand{\red}{\protect\bf}
\renewcommand\refname{\centerline
{\red {\normalsize \bf References}}}
\newtheorem{ass}
{\bf Assumption}[section]
\newtheorem{defi}
{\bf Definition}[section]
\newtheorem{tho}
{\bf Theorem}[section]
\newtheorem{rem}
{\sc Remark}[section]
\newtheorem{lemma}
{\bf Lemma}[section]
\newtheorem{coro}
{\bf Corollary}[section]
\newtheorem{prop}
{\bf Proposition}[section]
\renewcommand{\theequation}{\thesection .
\arabic{equation}}
\setcounter{section}{1}
\section*{\normalsize 1. Introduction}
\renewcommand{\theequation}{1.\arabic{equation}}
\renewcommand{\thetho}{1.\arabic{tho}}
\renewcommand{\thedefi}{1.\arabic{defi}}
\renewcommand{\therem}{1.\arabic{rem}}
\renewcommand{\theprop}{1.\arabic{prop}}
\renewcommand{\thelemma}{1.\arabic{lemma}}
\setcounter{equation}{0} \setcounter{lemma}{0} \setcounter{defi}{0}
\setcounter{prop}{0} \setcounter{rem}{0} \setcounter{tho}{0} This
note deals with the existence and regularity of solutions to  a singular
non-linear, second order, parabolic system, under Dirichlet
boundary conditions, of the type \be\label{PF}\begin{array}{ll}\dy
u_t-\nabla\cdot\left(|\nabla u|^{p-2}\nabla u\right)\dy= 0\,,&\hskip-0.2cm\textrm{
in }(0,T)\times\O,
\\\dy \hskip2.4cmu(t,x)\dy=0\,, &\hskip-0.2cm \textrm{ on }
(0,T)\times\po,\;\\\dy\hskip2.4cm u(0,x)=u_\circ(x),&\hskip-0.2cm\mbox{ on
}\{0\}\times\O,\end{array}\ee where
the $p$-growth exponent belongs
to the interval $(1,2)$. Here we
assume that $\O$ is a bounded
domain of $\R^n$, $n\geq2$, whose
boundary is $C^2$-smooth, and
$u:(0,T)\times\O\to \R^N,\,N\geq1$,
 is a scalar or a vector field.
The data $u_\circ$ belongs to $ L^2(\O)$.
 \par Our main purpose is
to prove ``global'', that is
on the cylinder $(\ve,T)\times\O$,
$\ve\geq0$, $L^r(\ve,T;L^q(\O))$ integrability
 properties of the second spatial
derivatives and of the time
derivative of solutions to problem
\eqref{PF}. Hence, for suitable
$p\in(1,2)$ and exponents $r,\,q$,
by Sobolev embedding theorems, we
deduce ``global'' regularity
of $u$ and $\nabla u$ in
H\"older spaces. Our results are
developed under two main
assumptions. The former is that we
 consider a $p$-parabolic system
with $p$-Laplacean operator and
not more general elliptic
operators, whose structural
properties have the $p$-Laplacean as
prototype. The latter concerns the
bounded domain. These assumptions
are made just  to develop in a
simpler way a new technique that
leads to the high integrability of
second derivatives and of the first
time derivate, which are the chief
results of this paper, and, as far
as we know, they are new in the
literature. The proof is performed
under the assumption of an
homogeneous right-hand side. This
is not a limit of the technique
employed, which works for a
nonnull right-hand side as well,
but a choice to develop the proofs in a more
readable way. As a
consequence of the integrability
properties, by embedding theorems,
we deduce the H\"older regularity
of the solutions to problem
\rf{PF}. This important topic has
been developed in a wide
literature. We refer  to the
 monograph \cite{DB} and
to the more recent essay
\cite{DBUV} for a general review,
and, for regularity results,
more specifically  to the papers
\cite{DBF,DBH,DBKV}.  For solutions to more
general singular or
degenerate parabolic
systems, as far as the local
integrability properties are concerned,
 we quote the papers \cite{AcMin,KL, LVP, Sch1, Sch2},
and, as far as the H\o lder regularity
is concerned, we quote the paper \cite{AMS} and
the very complete and recent memoirs
\cite{DMSt}.
\par Before giving the statement of
our results, we would like to say a
few words about the technique.
Firstly we point out that we do not
prove that a weak solution has more
regularity properties, but we prove
the existence of regular solutions,
and, as a consequence of the
uniqueness of weak solutions, the
same regularity is enjoyed by the
weak solution too. The
 existence theorem
is proved by using the Galerkin
method in the way suggested by
Prodi in \cite{prodi}, where the
special basis of eigenfunctions is
proposed, in our case of the
Laplacean operator. In this way,
provided that the initial data
$u_\circ\in L^2(\O)$, we are able
to furnish a solution to problem
\rf{PF} which has more regularity
of the usual weak solutions.
Indeed, we obtain $D^2u(t,x)\in
L^2(\ve,T; L^p(\O))$, $\ve\geq 0$.
This gives an advantage in
establishing further regularity
properties of the solutions. As far
as the integrability properties are
concerned, we employ a duality
technique which is a suitable
modification of the one employed in
\cite{MS} to estimate the second
derivatives in (space-time)
anisotropic
 Sobolev spaces. For
this task we define suitable
adjoint problems of \rf{PF}. This
approach firstly gives estimates
for the time derivatives $u_t$,
subsequently, as in \cite{MS},
viewing the first equation of
\rf{PF} as an equation of elliptic
type with data $f=u_t$, allows us
to establish estimates of the
second derivatives.
\par To better explain our results, we shortly
introduce the space
$V=W^{1,p}_0\cap L^2(\O)$. $V$ is a
reflexive Banach space endowed with
the norm
$\|\cdot\|_V=\|\cdot\|_{1,p}+\|\cdot\|_2$.
Moreover we denote by
$V'=W^{-1,p'}(\O)+L^2(\O)$ its
strong dual. For further notations
and positions see the next section.
  \begin{defi}\label{defnomu}
{\rm Let  $u_{\circ}\in L^2(\O)$. A field
$u\!:(0,T)\times \O\to\R^N$
 is said a solution of
system {\rm \eqref{PF}} if \be\label{ago151} u\in L^{p}(0,T; V)\cap
C(0,T;L^2(\O))\,, \ t^\frac 1p \nabla u\in
L^{\infty}(0,T;L^p(\O))\,,\ee \be\label{ago252}u_t\in
L^{p'}(0,T;V')\,,\ t\,u_t\in L^{\infty}(0,T;L^2(\O))\,,\
t^\frac{p+2}{2p}\,\nabla u_t \in L^{2}(0,T; L^p(\O))\,,\ee
$$\ba{ll}\dy \vs1\int_0^t [(u,\psi_\tau)-\left(|\nabla u|^{p-2}\,
\nabla u,\nabla \psi\right)]\,d\tau=(u(t),
\psi(t))-(u_\circ, \psi(0)),\\ \hskip6cm  \forall \psi\in
W^{1,2}(0,T;L^2(\O))\cap L^p(0,T;V),\ea$$ and
$$\lim_{t\to 0^+}\|u(t)-u_\circ\|_2=0\,.$$}
\end{defi}\par
The above definition is
different from the usual
formulation of a weak solution to
 problem \eqref{PF}, actually
the properties indicated for a
solution are wider than the ones
given in \cite{lions} and
considered by other authors, for instance
 in
\cite{DB,DMSt}. This is in
connection with the fact that we
are able to prove that, for all
$p\in(1,2)$, the set of solutions
is not empty,
 provided that $u_\circ\in
L^2(\O)$. Of course, a solution in
the sense of the
Definition\,\ref{defnomu} is a weak
solution in the sense given in
\cite{lions}. Hence the uniqueness
of the weak solution makes unique
the functional class of existence
and the related properties of the
solutions.\par We set
\be\label{restrp}p_\circ:=\max\big\{\mbox{\large$\frac
32, \frac{2n}{n+2}$}\big\}\,,\ee
and \be\label{P0} \ov
p:=2-\mbox{\large$\frac{1}{H}$}\,,\ee
with the constant $H$ introduced in
\eqref{C3rich}.
\begin{tho}\label{existencenomu}
{\sl Assume that $u_\circ$ belongs to
$L^2(\O)$. Then there exists a unique solution $u$ to problem
\eqref{PF} in the sense of Definition\,\ref{defnomu}.
Moreover, if $p>p_\circ$, then 
\be\label{sing1} t^{\alpha_1}\, \nabla u\in C(0,T;L^{2}(\O))\,,\ee
 and
\be\label{sing2}t^{\alpha_2}\, u\in L^{2}(0,T;W^{2,p}(\O))\,,\ee
where $\alpha_1=\frac{2}{2p-n(2-p)}$ and
$\alpha_2=\alpha_1+\frac{2-p}{2p}$.}
\end{tho}
\begin{coro}\label{corexistencenomu}
{\sl Assume that $u_\circ$ belongs to
$W_0^{1,2}(\O)$. Then there exists a unique solution $u$ to problem
\eqref{PF} in the sense of Definition\,\ref{defnomu}. Moreover, for
$p>\frac 32$,
$$\nabla u\in C(0,T; L^2(\O))\,,$$
and
$$ u\in L^{2}(0,T;W^{2,p}(\O))\,.$$}
\end{coro}
Theorem\,\ref{existencenomu} says
in particular
 that if the initial data is just in $L^2(\O)$, then
$\nabla u$, $D^2u$, $u_t$ and $\nabla u_t$ have
a singularity in the origin $t=0$, that we
explicitly compute. If the data is more regular,
as in Corollary\,\ref{corexistencenomu}, we
remove the singularity in $t=0$. This result
completely agrees with the known results for the
linear case, with an obvious rescaling due to
the exponent $p$.  In the case of a more regular
initial data, we limit ourselves to the
 claims in Corollary\,\ref{corexistencenomu} for the
sake of brevity. However under the assumption
$u_\circ \in W_0^{1,2}(\O)$, we could give
further regularity properties, that we consider
unessential for the developement of the paper.
We point out that these cannot be considered
like results on the asymptotic behavior of the
solution, since, as it is well known, for all
$p\!\in\! (1,2)$, if $\O$ is bounded there is
the extinction of the solution in a finite time
(cf. \cite{DB}).
\par We also observe that the introduction of a force term
 $f\in L^{p'}(0,T; V')$ on the right-hand side
  would be easy to handle and would lead to
 the same $L^2(\ve,T; W^{2,p}(\O))$, $\ve>0$, integrability for second
 derivatives. Obviously, under this weak assumption on $f$, the
 solution as in Definition\,\ref{defnomu} would lost the regularity
 properties of $u_t$ given in
 \eqref{ago252}$_{1,2}$.  \par
The next theorem and its corollary
are our chief results and concern
the ``global'' high regularity of
the solutions furnished by
Theorem\,\ref{existencenomu}. For
the definition of the H\"{o}lder
seminorm in \eqref{holdertx} we
refer to the next section. We set
$$p_1:=
\mbox{\large$ \frac{7(n-2)+1-
\sqrt{4(n-1)^2-3}}{3(n-2)}$}\,.$$
\begin{tho}\label{regularitynomu}
{\sl Let $n\geq 3$. Let $p>
\max\{p_\circ,p_1\}$, and let $q\in
[\frac{2n}{n(p-1)+2(2-p)},
\frac{7-3p}{4-2p}]$. Assume that
$\O$ is a convex domain. If
$u(t,x)$ is the solution of Theorem\,\ref{existencenomu}, then, for all
$\ve>0$,\, \be\label{CRE} u\in
L^\infty(\ve,T;W^{2,\widehat
q}(\O))\mbox{ with }u_t\in
L^\infty(\ve,T;L^q(\O)),\ee with
$\widehat
q=\frac{nq(p-1)}{n-q(2-p)}$ if
$q<n$, $\widehat q<n$ if $q=n$, and
$\widehat q=q$ if $q>n$.
\par Under the further
assumption $p>\ov p$, the same
result holds for $\O$ non-convex
domain.}
\end{tho}
\begin{tho}\label{regularitynomun2}
{\sl Let $n=2$. Let $p>\frac 32$,
and let $q\in (2,
\frac{7-3p}{4-2p}]$. Assume that
$\O$ is a convex domain. If
$u(t,x)$ is the solution of
Theorem\,\ref{existencenomu}, then, for all
$\ve>0$,\, \be\label{CREn2} u\in
L^\infty(\ve,T;W^{2, q}(\O))\mbox{
with }u_t\in
L^\infty(\ve,T;L^q(\O)).\ee
\par Under the further
assumption $p>\ov p$, the same
result holds for $\O$ non-convex
domain.}
\end{tho}We set
$$p_2:=\mbox{\large$\frac{2n+7-
\sqrt{(2n-7)^2+8n}}6$}.$$
\begin{coro}\label{corregularitynomu}
{\sl Assume that $\O$ is a convex
domain and $u(t,x)$ is the solution
of Theorem\,\ref{existencenomu}. Let
$p>\frac32$ for $n=3$ and $p>p_2$
for $n>3$, and $q_0\in(\frac
np,n]$.  Then, for each
$t>t_\circ>0$ we get
\be\label{holdertxa}
\big[\,u\,\big]_{\lambda_\circ,t,x}\leq c\,(
t_\circ ^{-1-\gamma_\circ}
\|u_\circ\|_2^{(2-p)\gamma_\circ+1}+
t_\circ
^{-\frac{1+\gamma_\circ}{p-1}}
\|u_\circ\|_2^{\frac{(2-p)\gamma_\circ+1}{p-1}}
)\,,\ee where $\lambda_\circ=2-\frac
{n}{\widehat q_\circ}$, with $\widehat
q_\circ=\frac{nq_\circ(p-1)}{n-q_\circ(2-p)}$
if $q_\circ<n$, $\widehat
q_\circ\in (\frac n2, n)$ if
$q_\circ=n$, and
$\gamma_\circ=\frac{n(q_\circ-2)}{q_\circ(
 2p-2n+np)}$. Moreover, let
 $p>\max\{p_\circ,
 \frac{4n-7}{2n-3}\}$ and $q_1\in (n,
 \frac{7-3p}{4-2p}]$. Then, we get
 \be\label{holdertx}
\big[\,\nabla u\,\big]_{\lambda_1,t,x}\leq c\,( t_\circ
^{-1-\gamma_1}
\|u_\circ\|_2^{(2-p)\gamma_1+1}+
t_\circ ^{-\frac{1+\gamma_1}{p-1}}
\|u_\circ\|_2^{\frac{(2-p)\gamma_1+1}{p-1}}
)\,.\ee
where $\lambda_1=1-\frac
{n}{q_1}$,  and
$\gamma_1=\frac{n(q_1-2)}{q_1(2p-2n+np)}$.
The constant $c$ in
\rf{holdertxa}-\rf{holdertx} is
independent of $t_\circ$ and
$u_\circ$.\par Under the further
assumption $p>\ov p$, the same
results hold for $\O$ non-convex
domain.}
\end{coro}
\begin{coro}\label{corregularitynomu2}
{\sl Let $n=2$. Let $p>\frac 32$,
and $q\in(2, \frac{7-3p}{4-2p}]$.
Assume that $\O$ is a convex domain
and $u(t,x)$ is the solution of
Theorem\,\ref{existencenomu}. Then,
for each $t>t_\circ>0$ we get
 \be\label{holdertxb}
\big[\,\nabla u\,\big]_{\lambda,t,x}\leq c\,( t_\circ
^{-1-\gamma_1}
\|u_\circ\|_2^{(2-p)\gamma_1+1}+
t_\circ ^{-\frac{1+\gamma_1}{p-1}}
\|u_\circ\|_2^{\frac{(2-p)\gamma_1+1}{p-1}}
)\,.\ee
where $\lambda=  1-\frac
{2}{q}$, and
$\gamma_1=\frac{(q-2)}{2q(p-1)}$.
The constant $c$ in \rf{holdertxb}
is independent of $t_\circ$ and
$u_\circ$.\par Under the further
assumption $p>\ov p$, the same
results hold for $\O$ non-convex
domain.}
\end{coro}
Our result of ``high regularity''
is expressed by means of the
existence of the second
derivatives, that is estimates
\eqref{CRE} and \eqref{CREn2} .
These estimates for suitable $p$
and $q$  imply the H\"older
regularity of $u$ and $\nabla u$.
We point out that our H\"older
exponent $\lambda$ depends on $p,
n$. This is in accordance with the
result given in \cite{CDB}. However
we are not able to compare the two
exponents since, as far as we
know, in \cite{CDB}
a functional dependence  for $\lambda$  on $p$ is
not given.
\par We have recently seen paper \cite{BDV},
where, for $n\geq 3$, under the
assumption of $f\in
L^{p'}(0,T;W^{-1,p'}(\O))\cap
L^2(0,T;L^2(\O))$, it is proved
that a weak solution of \eqref{PF}
belongs to $L^{2(p-1)}(0,T;
W^{2,r}(\O))$, with
$r=\frac{2n(p-1)}{n-2(2-p)}$ under
suitable constraints on $p,r$.
Hence, for $n=3$,  $r$ belongs to
$(p,2)$,  and, for $n>3, r<p$.
\par Finally, we prove a global
pointwise bound for the solution.
\begin{tho}\label{mmtp}{\sl  Let $u$ be the solution of \eqref{PF} corresponding to
$u_\circ\in L^\infty(\O)$. Then
\be\label{mmtp1}\|u(t)\|_\infty\leq \|u_\circ\|_\infty.\ee Moreover,
if $p>\frac{2n}{n+2}$ then, corresponding to an
initial data  $u_\circ\in L^q(\O)$, for some $q\in [2,+\infty]$ one
has \be\label{mmtp2}\|u(t)\|_\infty\leq c\,\|u_\circ\|_q\,
\|u_\circ\|_2^\frac{2(2-p)\beta}{q}  t^{-\frac{2\beta}{q}}  , \quad
\forall t>0,\ee with $\beta:=\frac{n}{p(n+2)-2n}$.
 }
\end{tho}
The precise aim of Theorem\,\ref{mmtp} is to prove a
$L^\infty(\O)$-bound for a weak
solution with no investigations of
high regularity properties of
solutions. Of course, estimate
\eqref{mmtp2} holds for $t\in
(0,T)$, where $T$ is the instant of
extinction. Analogous results are
proved in \cite{DBH, DBUV,DB} for
equations, and locally in
\cite{choe} for systems.\vskip0.1cm
Finally, we shortly describe the
plan of the paper and the strategy
of the proofs. In Sec.\,2, we
introduce some notations and
auxiliary lemmas. In particular we
prove Lemma\,\ref{LL1}, which is an
important tool to estimate the
second derivatives, and semigroup
properties for a suitable linear
parabolic system with regular
coefficients. In Sec.\,3 we
introduce two approximating
systems, \eqref{PFl} and
\eqref{PFepv}. They are both
non-singular quasi-linear systems,
for which the Galerkin
approximation method, and in
particular a suitable choice of the
basis for this approximation,
together with weighted estimates in
$W^{m,r}(\O)$ and Lemma\,\ref{LL1} are
the essential tools to get
existence and regularity. Since the
proof of the basic properties of
these systems is standard for
people acquainted with the Galerkin
method, we confine it in the
appendix, at the end of the paper.
The proof of our existence result
is then given in Sec.\,4. The
crucial step in the proof of the
higher integrability of second
derivatives is the derivation of an
$L^\infty(0,T; L^q(\O))$ estimate
on the time derivative $u_t$, which
is done in Sec.\,5, by using the
semigroup properties of Sec.\,2.
Once this regularity has been
derived, the
 corresponding $L^\infty(0,T; W^{2,\widehat q}(\O))$-integrability of
 the second derivatives is obtained in Sec\,6, and it relies on
 the regularity results on the $p$-Laplacean elliptic system
studied in
 \cite{CMellittico}. For $\widehat q>n$ this result gives us the H\"{o}lder continuity of the solution.
  Finally, the maximum modulus theorem is proved in Sec.\,7,
  employing a duality arguments.
\section*{\normalsize 2. Notations and some auxiliary results}
\renewcommand{\theequation}{2.\arabic{equation}}
\renewcommand{\thetho}{2.\arabic{tho}}
\renewcommand{\thedefi}{2.\arabic{defi}}
\renewcommand{\therem}{2.\arabic{rem}}
\renewcommand{\theprop}{2.\arabic{prop}}
\renewcommand{\thelemma}{2.\arabic{lemma}}
\renewcommand{\thecoro}{2.\arabic{coro}}
\setcounter{equation}{0} \setcounter{coro}{0} \setcounter{lemma}{0}
\setcounter{defi}{0} \setcounter{prop}{0} \setcounter{rem}{0}
\setcounter{tho}{0}
Throughout the paper we denote by $p$ the growth exponent, with $p\in (1,2)$.
We denote by $\O\subset\R^n$ a bounded domain
whose boundary $\po$ is $C^2$-smooth. For a function $v(t,x)$,
by $\pa_k v$ and $\pa_t v$ we mean $\frac{\pa}{\pa x_k}v(t,x)$
 and  $\frac{\pa}{\pa t}v(t,x)$, respectively.
 We set $v \cdot\nabla v=v_k\partial_k\, v$, and
$\nabla v \cdot v=v\cdot
\partial_k\, v$. For $m\in \N\cup\{0\}$,
$C^m(\ov\O)$ ($C(\ov\O)$ for $
m=0$) is the usual space of
functions which are bounded and
uniformly continuous on $\O$
together with their derivatives up
to the order $m$. The norm in
$C^m(\ov\O)$ is denoted by $|\cdot
|_m:={\overset
m{\underset{|\alpha|=0}\sum}}{\underset\O\sup}|D^\alpha
u(x)|$. For $\lambda\in(0,1)$, by
$C^{m,\lambda}(\ov\O)$ we mean the
set of functions of $C^m(\ov\O)$
such that, for $|\alpha|=m$,
$D^\alpha u\in
C^{0,\lambda}(\ov\O)$, that is
$|D^\alpha u|_{0}+[D^\alpha
u]_{\lambda}<\infty$, where
$[\,\cdot\,]_{\lambda}$ is the
H\"older seminorm. The norm of an
element of $ C^{m,\lambda}(\ov\O)$
is denoted by
$|u|_{m,\lambda,\O}:=|u|_m+[D^\alpha
u]_{\lambda}\,, |\alpha|=m$. We
denote by $C^m(a,b;X)$ the Banach
space (endowed with the natural
norm) of all functions bounded and
continuous on $(a,b)\subseteq\R$
with value in a Banach space $X$,
together with all derivatives
$D^k$, $k\leq m$. For $\lambda\in
(0,1)$ we set \be\label{posalpha}
[\,g\,]_{\lambda,t,x} =\sup_{\ov x,
\ov {\ov x}\in \O, \  \ov
x\not=\ov{\ov x} \atop t\in (0,T)}
 \frac{|g(t, \ov x)-g(t,\ov {\ov x})|}{|\ov x-\ov {\ov x}|^\lambda}+
 \sup_{x\in \O, \atop
 \ov t, \ov{\ov t}\in (0,T),\  \ov t\not=\ov{\ov t}}
 \frac{|g(\ov t, x)-g(\ov {\ov t}, x)|}{|\ov t-\ov {\ov t}|^\frac\lambda 2}\,,\ee
provided that the right-hand side is finite. The $L^p$-norm is
denoted by $\|\cdot\|_p$ and, if $m\geq0$, the $W^{m,p}$ and
$W_0^{m,p}$-norms are denoted by $\|\cdot\|_{m,p}$. We introduce the
space $V=W^{1,p}_0\cap L^2(\O)$. $V$ is a reflexive Banach space
endowed with the norm $\|\cdot\|_V=\|\cdot\|_{1,p}+\|\cdot\|_2$,
where $\|\cdot\|_{1,p}$ represents a semi-norm on $V$. Moreover we
denote by $V'=W^{-1,p'}(\O)+L^2(\O)$ its strong dual. Note that
$W^{1,p}_0\subset L^2(\O)$ only if $p\geq \frac{2n}{n+2}$. On the
other hand $V$ is dense and continuously embedded in $L^2$ and
$V\subset L^2\subset V'$. Let $q\in[1,\infty)$, let $X$ be a Banach
space with norm $\|\cdot\|_X$. We denote by $L^q(a,b;X)$ the set of all
function $f:(a,b)\to X$ which are measurable and such that the
Lebesgue integral \mbox{\footnotesize $\dy\intll ab$}$\|
f(\tau)\|^q_Xd\tau= \|f\|_{L^q(a,b;X)}<\infty$. As well as, if
$q=\infty$ we denote by $L^\infty(a,b;X)$ the set of all function
$f:(a,b)\to X$ which are
 measurable and such that ${\rm ess\ sup}_{t\in(a,b)}\ \|
f(t)\|^q_X=\|f\|_{L^\infty(a,b;X)}<\infty$.  
\par
In the remaining part of this section we give some preliminary
results,
 which represent fundamental tools in our proofs. The first
is  the following lemma, which, for $p=2$, gives a well known
estimate (see \cite{L} and \cite{LU}).
\begin{lemma}\label{LL1}{\sl
Let $\mu>0$. Assume that $v \in
W^{2,2}(\O)\cap W_0^{1,2}(\O)$. Then, for any $\eta>0$,
$$\big\|(\mu+|\nabla v|^2)^\frac{(p-2)}{4} {D^2 v}\big\|_2\,
\leq C_1\big\|(\mu+|\nabla v|^2)^\frac{(p-2)}{4}{\Delta v}
\big\|_2\!\!+\frac{C_2}{\eta}\left(\|\nabla v\|_p^p+\mu^\frac p2
|\Omega|\right)^\frac 12.$$ where
$$C_1:=\mbox{\large{$\left(\frac{p}{p(p-1)^2-\eta}\right)^{\!\frac12}$}}.$$
If $\O$ is a convex domain the inequality holds with $C_2=0$,
$\eta=0$.}
 \end{lemma}
 \Pr  We prove the result for sufficiently smooth functions. It can be extended to functions
 in $W^{2,2}(\O)\cap W_0^{1,2}(\O)$ by density arguments.
 So, let $v$ be a function which is continuously differentiable three
times and vanishes on $\po$. Integration by parts gives
$$\ba{ll}\vs1 \dy \int_\O (\mu+|\nabla
v|^2)^\frac{(p-2)}{2} |\Delta v|^2\, dx =-\int_\O (\mu+|\nabla
v|^2)^\frac{(p-2)}{2}\frac{
\partial\Delta v}{\partial x_k}\cdot\frac{\partial v}{\partial
x_k}\,dx\\\hfill \vs1 \dy -(p-2)\int_\O (\mu+|\nabla
v|^2)^\frac{(p-4)}{2}\Delta v\cdot\frac{\partial v}{\partial x_k}\,
\nabla v\cdot\frac{\partial\nabla v}{\partial x_k}\,dx \dy
+\int_{\po} (\mu+|\nabla v|^2)^\frac{(p-2)}{2}\Delta
v\cdot\frac{\partial v}{\partial n}\,d\sigma\\\hfill\vs1 \dy
=\int_\O\!(\mu+|\nabla v|^2)^\frac{(p-2)}{2}\!\frac{\partial
v}{\partial x_j\partial x_k}\cdot\frac{\partial v}{\partial
x_j\partial x_k}\,dx + (p-2)\!\int_\O\! (\mu+|\nabla
v|^2)^\frac{(p-4)}{2}\!\left(\frac{\partial \nabla v}{\partial
x_j}\cdot \nabla v\!\right)^2\!dx\\\vs1 \dy\hskip0.5cm
 -(p-2)\int_\O (\mu+|\nabla
v|^2)^\frac{(p-4)}{2}\Delta v\cdot\frac{\partial v}{\pa x_k}\,\nabla
v\cdot\frac{\pa\nabla v}{\pa x_k}\,dx \dy
\\\hskip0.5cm \dy+\int_{\po} (\mu+|\nabla
v|^2)^\frac{(p-2)}{2}\left[\Delta v\cdot\frac{\partial v}{\partial
n}-\frac{\partial v}{\partial x_k\partial n}\cdot\frac{\partial
v}{\pa x_k} \right] \,d\sigma \,.\ea$$ Denote the boundary integral
in the previous estimate by $I_{\po}$. Since $p>1$, one can estimate the right-hand side
as follows \be\label{L2a}\ba{ll}\vs1 \dy (p-1)\int _\O (\mu+|\nabla
v|^2)^\frac{(p-2)}{2} |D^2 v|^2\, dx\\\dy \leq\!\int_\O (\mu+|\nabla
v|^2)^\frac{(p-2)}{2}|\Delta v|^2\,dx+ (2-p)\!\int_\O (\mu+|\nabla
v|^2)^\frac{(p-2)}{2} |D^2 v|\, |\Delta v|\,dx- I_{\po}\,.\ea\ee

\par {\it $\O$ convex} - By
using the arguments in \cite{L}, based on a localization technique,
one can show that the boundary integral $I_{\po}$ is non-negative if
$\O$ is convex. Therefore from \eqref{L2a} one gets $$\ba{ll}\vs1
\dy (p-1)\int _\O (\mu+|\nabla v|^2)^\frac{(p-2)}{2} |D^2 v|^2\,
dx\\\dy \leq \int_\O (\mu+|\nabla v|^2)^\frac{(p-2)}{2}|\Delta
v|^2\,dx+ (2-p)\int_\O (\mu+|\nabla v|^2)^\frac{(p-2)}{2} |D^2 v|\,
|\Delta v|\,dx\,.\ea$$ By applying H\"{o}lder's and  Cauchy's
inequalities to the last integral one readily has
 \be\label{L3a}\ba{ll}\dy
\left[p-1-\frac{\ve}{2}(2-p)^2\right]\!\int _\O (\mu+|\nabla
v|^2)^\frac{(p-2)}{2} |D^2 v|^2\, dx \\ \dy \hskip4cm \leq
\left(\!1+\frac{1}{2\ve}\right)\! \int_\O (\mu+|\nabla
v|^2)^\frac{(p-2)}{2}|\Delta v|^2\,dx\,,\ea \ee hence \be\label{L3} \int
_\O (\mu+|\nabla v|^2)^\frac{(p-2)}{2} |D^2 v|^2\, dx \leq C(\ve)
\int_\O (\mu+|\nabla v|^2)^\frac{(p-2)}{2}|\Delta v|^2\,dx\,,\ee
with $$C(\ve):=\mbox{\large${\frac{1+2\ve}{\ve[2(p-1)-\ve(2-p)^2]}}$}\,.$$
By an
easy computation, one can verify that the minimum of $C(\ve)$ equals
$1/(p-1)^2$ and it is attained for $\ve=(p-1)/(2-p)$. Therefore we
get \be\label{L4} \int _\O (\mu+|\nabla v|^2)^\frac{(p-2)}{2} |D^2
v|^2\, dx \leq \frac{1}{(p-1)^2} \int_\O (\mu+|\nabla
v|^2)^\frac{(p-2)}{2}|\Delta v|^2\,dx\,.\ee  \par {\it $\O$
non-convex } - If $\O$ is not convex, starting from \eqref{L2a} and
using  the above arguments (see \eqref{L3a}--\eqref{L4}) we have
\be\label{L4a} \int _\O (\mu+|\nabla v|^2)^\frac{(p-2)}{2} |D^2
v|^2\, dx \leq \frac{1}{(p-1)^2} \int_\O (\mu+|\nabla
v|^2)^\frac{(p-2)}{2}|\Delta v|^2\,dx-
\frac{2}{p(p-1)}\,I_{\po}\,.\ee Again following \cite{L}, the
integral $I_{\po}$ can be estimated as follows
\be\label{L5}\ba{ll}\dy \vs1 I_{\po}\!\!&\dy \leq
C\int_{\po}(\mu+|\nabla v|^2)^\frac{(p-2)}{2} \left(\frac{\partial
v}{\partial n}\right)^2 d\sigma\leq C \|\nabla v\|_{L^p(\po)}^p\\
&\dy  \leq C\int_\O|\nabla v|^pdx + C\int_\O|\nabla
v|^{p-1}|D^2v|\,dx\,.\ea\ee Multiplying and dividing by
$(\mu+|\nabla v|^2)^\frac{(2-p)}{4}$, using H\"{o}lder's and then
Cauchy's inequalities  we have, for any $\eta>0$,
$$\ba{ll}\dy \vs1 \int_\O\!|\nabla v|^{p-1}|D^2v|\,dx\!\!\!\!& \dy\leq
\left(\int_\O\! (\mu+|\nabla v|^2)^\frac{2-p}{2} |\nabla
v|^{2(p-1)}dx\!\right)^{\!\frac 12}\!\!\left(\int_\O\! (\mu+|\nabla
v|^2)^\frac{p-2}{2}|D^2v|^2dx\!\right)^{\!\frac 12}\\&\dy  \leq
\frac{1}{2\eta} \int_\O (\mu+|\nabla v|^2)^\frac{p}{2} dx+
\frac{\eta}{2}\int_\O (\mu+|\nabla v|^2)^\frac{p-2}{2}|D^2v|^2dx
\,,\ea$$ for any $\eta>0$. Therefore,
$$I_{\po}\leq C \left(\|\nabla v\|_p^p+\frac{1}{2\eta}\| (\mu+|\nabla v|^2)^\frac 12\|_p^p
 + \frac{\eta}{2} \| (\mu+|\nabla v|^2)^\frac{(p-2)}{4}
D^2v\|_2^2\right) \,.
$$
By replacing the above estimate in  \eqref{L4a} we get
$$\ba{ll}\vs1 \dy \int _\O (\mu+|\nabla v|^2)^\frac{(p-2)}{2} |D^2
v|^2\, dx \leq \frac{1}{(p-1)^2}\int_\O (\mu+|\nabla
v|^2)^\frac{(p-2)}{2}|\Delta v|^2\,dx\\\hskip0.5cm \dy
+\frac{\eta}{p(p-1)}\int _\O (\mu+|\nabla v|^2)^\frac{(p-2)}{2} |D^2
v|^2\, dx +C\,\left( \|\nabla v\|_p^p+\frac{1}{2\eta}\| (\mu+|\nabla
v|^2)^\frac 12\|_p^p\right),\ea$$ which easily gives the
result.\chiu \vskip0.5cm Our second kind of results is concerned
with the analysis of semigroup properties for the following
parabolic system with regular coefficients
\be\label{AD1}\begin{array}{ll}\dy\vs1 \vp_s-\nu\Delta \vp- \nabla\cdot
(B_\eta(s,x)\nabla \vp)= 0\,,&\hskip-0.2cm\textrm{ in }(0,t)\times\O,
\\\dy\vs1\hskip3.5cm \vp(s,x)=0\,,&\hskip-0.2cm\textrm{ on }
(0,t)\times\po,\;\\\dy\hskip3.5cm \vp(0,x)=\vp_\circ(x),&\hskip-0.2cm\mbox{ on
}\{0\}\times\O,\end{array}\ee with $\nu\geq 0$ and
$B_\eta(s,x)=(B_\eta)_{i\alpha j\beta}(s,x)$ satisfying the
following conditions \be\label{coerc}\ba{ll}\vs1\dy  B_\eta\ \mbox
{is continuous in }  [0,t]\times \overline\O, \\\dy
\|B_\eta\|_{\infty}\equiv\max_{i,\alpha,j,\beta}\|(B_\eta)_{i\alpha
j\beta}\|_{\infty}<+\infty,\\
B_\eta(s,x) \mbox{ is uniformly elliptic}\,. \ea\ee
\begin{lemma}\label{newag1}{\sl Assume that $\nu>0$ and let $\vp_\circ(x)\in C_0^{\infty}(\O)$.  Then, there exists
a unique solution $\vp$ of \eqref{AD1}, such that $\vp\in
 L^2(0,t; W^{2,2}(\O)\cap W_0^{1,2}(\O))$, $\vp_s\in
L^2(0, t; L^{2}(\O))$.
 }\end{lemma}
\Pr The existence and regularity follow from well known regularity
results for linear parabolic systems with uniformly continuous and
bounded coefficients. We refer, for instance, to \cite{LSU}, Theorem
IV.9.1.\chiu

\vskip0.2cm
For  $p\in (1,2)$ and $\mu>0$, set
\be\label{oper3}
a(\mu,v):=\left(\mu+|(\nabla
v)|^2\right)^\frac{(p-2)}{2}\,,\ee
and
\be\label{amu}
a_\eta(\mu,v):=\left(\mu+|J_\eta(\nabla
v)|^2\right)^\frac{(p-2)}{2}\,,\ee
 with
$J_{\eta}$ space-time Friederich's mollifier,
and assume
  that
 \be\label{beeta2}
 s \,\|(\mu+|\nabla v|^2)^\frac 12\|^p_{p} \leq
 M\,,\quad \forall s\in [0,t],\ee
with a positive constant $M$.\par
Let   $p\in (1,2)$,  $\mu>0$, and define
\be\label{beeta}\ba{ll}\dy\vs1  (B_\eta(s,x))_{i\alpha j\beta}:=
\frac{\delta_{ij}\,\delta_{\alpha\beta}}{(\mu+|J_\eta(\nabla
v)(t-s,x)|^2)}{\null_{\frac {2-p}{2}}}\\
\dy\hskip5cm-b\,(2-p)\frac{(J_\eta(\nabla v)\otimes J_\eta(\nabla
v))(t-s,x)}{(\mu+|J_\eta(\nabla v)(t-s,x)|^2)^{\frac
{4-p}{2}}}\ea\ee with
  $J_{\eta}$ space-time Friederich's mollifier and $b=0,1$.
  Note that $B_\eta$ defined in \eqref{beeta}, by means of \eqref{amu}, satisfies condition \eqref{coerc}, with
$$\|B_\eta\|_{\infty}<(3-p) \mu^\frac{p-2}{2}<+\infty,\ \forall
\mu>0\,.$$

\begin{lemma}\label{adgt}{\sl
 Assume that $\vp_\circ(x)\in C_0^{\infty}(\O)$
  and let $\vp$ be the unique solution $\vp$ of \eqref{AD1}, corresponding to $B_\eta$
  as in \eqref{beeta}. Then, for all $r\in [1,2]$ if $b=0$,
  and for all $r\in [\frac{7-3p}{3-p},2]$ if $b=1$,
  \be\label{max11a}\|\vp(s)\|_r
  \leq \|\vp_\circ\|_r,\ \forall s\in [0,t], \ \textrm{uniformly in } \nu>0 \textrm{ and } \eta>0\,.
\ee
 }\end{lemma}
 \Pr
Let us multiply \eqref{AD1}$_1$ by $\vp (\delta+|\vp|^2)^\frac
{r-2}{2}$, for some $\delta>0$. Then \be\label{max13}\ba{ll}\vs1 \dy
\frac 1r\frac{d}{ds}\|(\delta+|\vp|^2)^\frac 12\|_r^r+\nu \int_\O
(\delta+|\vp|^2)^\frac {r-2}{2}|\nabla\vp|^2dx\\ \hskip0.5cm \vs1\dy
+\nu(r-2)\!\!\int_\O (\delta+|\vp|^2)^\frac {r-4}{2}(\nabla\vp\cdot
\vp)^2dx \dy + \int_\O a_\eta(\mu,v)\,(\delta+|\vp|^2)^\frac
{r-2}{2}\,|\nabla\vp|^2dx
\\\dy\hskip0.5cm \vs1  +(r-2)\!\! \int_\O a_\eta(\mu,v)\,(\delta+|\vp|^2)^\frac {r-4}{2}\,(\nabla\vp\cdot \vp)^2 dx\\
 \vs1\hskip0.5cm \dy+ b\, (p-2) \bigg[ \int_\O (\mu+|J_\eta(\nabla
v)|^2)^\frac{p-4}{2}(\delta+| \vp|^2)^\frac {r-2}{2} (J_\eta(\nabla
v)\cdot \nabla \vp)^2\,
dx \\
\dy \hskip0.5cm+ (r-\!2)\!\!\int_\O \frac{(\delta+|\vp|^2)^\frac
{r-4}{2}}{(\mu+|J_\eta(\nabla v)|^2)}{\null_\frac{4-p}{2}}\, (J_\eta(\nabla v)\cdot\!
\nabla\vp) (J_\eta(\nabla
v)\cdot\vp)\, (\nabla \vp\cdot\vp)dx\bigg]=0.\ea\ee Taking into account that
$r\geq 1$, the differential equation \eqref{max13} gives $$\ba{ll}\vs1 \dy \frac
1r\frac{d}{ds}\|(\delta+|\vp|^2)^\frac 12\|_r^r\dy+ (r-1)\int_\O a_\eta(\mu,v)(\delta+|
\vp|^2)^\frac {r-2}{2}\,|\nabla \vp|^2dx\\\dy\vs1 \hskip2cm  \leq b\,
(2-p)(3-r)\int_\O a_\eta(\mu,v)(\delta+|\vp|^2)^\frac
{r-2}{2}\,|\nabla \vp|^2 dx \,.\ea$$ Therefore, for all $r\in [1,2]$ if $b=0$, and for all $r\in [
\frac{7-3p}{3-p},2]$ if $b=1$, a straightforward computation gives the existence of a constant $C\geq 0$ such that
$$\ba{ll}\vs1 \dy \frac 1r\frac{d}{ds}\|(\delta+|\vp|)^\frac 12\|_r^r
+ C \int_\O a_\eta(\mu,v)(\delta+|\vp|^2)^\frac {r-2}{2}\,|\nabla
\vp|^2dx \leq 0\,,\ea$$ from which one obtains $$
\|(\delta+|\vp(s)|^2)^\frac 12\|_r\leq
\|(\delta+|\vp_\circ|^2)^\frac 12\|_r, \quad \forall s\in [0,t], \
\forall \delta>0\,,$$ which gives \eqref{max11a}.\chiu
\begin{lemma}\label{newag2}{\sl Let $p>\frac{2n}{n+2}$. Under the assumptions of Lemma\,\ref{adgt}
\be\label{Maxpostaa}\|\vp(s)\|_2\leq c\,{M^{\frac{(2-p)\gamma}{2}}\,}\|\vp_\circ\|_r\,  [t^\frac
1p-(t-s)^\frac 1p]^{-p\gamma}\,,\quad \forall s\in\, (0,t], \ee
\be\label{Maxpostraa} \gamma=\gamma(r):=\mbox{\large$
\frac{n(2-r)}{2rp-2nr+rnp}$}\,.\ee Moreover, if $b=0$, then
\be\label{Maxpost}\|\vp(s)\|_r\leq  c\, M^\frac{(2-p)\beta}{r'} \|\vp_\circ\|_1\,
[t^\frac 1p-(t-s)^\frac
1p]^{-2p\beta\left(1-\frac{1}{r}\right)}\,,\quad \forall
s\in\,(0,t]\,,\ee with \be\label{Maxpostr} \beta :=\gamma(1)=\mbox{\large$\frac{n}{p(n+2)-2n}$}\,.\ee
 and $r'$ conjugate exponent of $r$.
   }\end{lemma}

 \Pr  Let us first observe that,
multiplying equation \eqref{AD1}$_1$ by $\vp$, taking into account
 \eqref{beeta} and integrating over $\O$, we
obtain \be\label{max4}\frac 12\frac{d}{ds}\|\vp\|_2^2+\nu \int_\O
|\nabla \vp|^2dx+ C\int_\O a_\eta(\mu,v)\,|\nabla \vp|^2dx\leq
0\,,\ee hence \be\label{max5}\|\vp(s)\|_2^2+2\nu\int_0^s\|\nabla
\vp(\tau)\|_2^2\,d\tau+ 2C\,\int_0^s \|a_\eta(\mu,
v(t-\tau))^\frac{1}{2}\,\nabla \vp(\tau)\|_2^2\,d\tau \leq
\|\vp_\circ\|_2^2\,.\ee By H\"{o}lder's inequality
\be\label{maxa1}\ba{ll}\dy\vs1 \int_\O |\nabla \vp|^pdx= \int_\O
(\mu+J_\eta|\nabla v|^2)^\frac{p(p-2)}{4}\,|\nabla
\vp|^p(\mu+|J_\eta(\nabla v)|^2)^\frac{p(2-p)}{4}dx\\\vs1 \dy \leq
\left(\int_\O (\mu+|J_\eta(\nabla v)|^2)^\frac{(p-2)}{2}\,|\nabla
\vp|^2\right)^\frac p2 \left(\int_\O (\mu+|J_\eta(\nabla
v)|^2)^\frac{p}{2}dx\right)^\frac{2-p}{2}\\\dy \leq
\|(\mu+|J_\eta(\nabla v)|^2)^\frac{(p-2)}{4}\,\nabla \vp\|_2^p\,
\|(\mu+|J_\eta(\nabla
v)|^2)\|_{\frac p2}^\frac{p(2-p)}{4}\,.\ea\ee%
By using Minkowski's inequality, the last term can be treated as
follows
$$\|(\mu+|J_\eta(\nabla
v)(t-s)|^2)\|_{\frac p2}^\frac 12\leq \|J_\eta(\sqrt {\mu}+|\nabla
v|)(t-s)\|_{p}\leq \!\int_\R\!\! J_\eta(t-s-\tau) \|\sqrt {\mu}+|\nabla
v(\tau)|\|_{p}\,d\tau.$$ Therefore, using \eqref{beeta2}, it can be
further estimated as
$$\|(\mu+|J_\eta(\nabla
v)(t-s)|^2)\|_{\frac p2}^\frac 12 \leq 2 \int_\R\! J_\eta(t-s-\tau)
\|(\mu+|\nabla v(\tau)|^2)^\frac 12\|_{p}\,d\tau \leq
 c\left(\frac{M}{t-s}\right)^{\!\frac 1p}\!,$$ which, raised to the
 power $(2-p)\,$ gives, a.e. in $s>0$,
$$\ba{ll}\vs1\dy
\|(\mu+|J_\eta(\nabla v)(t-s)|^2)\|_{\frac p2}^\frac{2-p}{2}\leq
c \left(\frac {M}{t-s}\right)^{\!\frac{\!2-p}{p}} \,.\ea$$ Using this
estimate in \eqref{maxa1} we end-up with
\be\label{maxxx4}c\left(\frac{t-s}{M}\right)^{\!\frac{2-p}{p}}\|\nabla
\vp\|_p^2\leq \|(\mu+J_\eta|\nabla v|^2)^\frac{(p-2)}{4}\,\nabla
\vp\|_2^2\,. \ee From the differential inequality \eqref{max4} and
\eqref{maxxx4} we obtain \be\label{maxa2aa}\frac
12\frac{d}{ds}\|\vp\|_2^2+ c\,
\left(\frac{t-s}{M}\right)^{\!\frac{2-p}{p}}\|\nabla \vp\|_p^2\leq
0\,.\ee By the well known Gagliardo-Nirenberg inequality and estimate
\eqref{max11a} we have \be\label{maxa3aa}\|\vp\|_2\leq c\|\nabla
\vp\|_p^a\,\|\vp\|_r^{1-a}\leq c\|\nabla
\vp\|_p^a\,\|\vp_\circ\|_r^{1-a}\,,\ee with $
a=\frac{np(2-r)}{2(np+rp-nr)}\,.$
  From \eqref{maxa2aa} and \eqref{maxa3aa} we arrive at the differential inequality
$$\frac 12\frac{d}{ds}\|\vp\|_2^2+ C\,
\left(\frac{t-s}{M}\right)^{\!\frac{2-p}{p}}\|\vp_\circ\|_r^\frac{2(a-1)}{a}
\|\vp\|_2^\frac 2a\leq 0\,\quad  \forall s\in[0,t].$$
 Integrating from $0$ to
$s$ and performing straightforward calculations, we get
$$\|\vp(s)\|_2\leq c M^{\frac{(2-p)}{2p}\frac{a}{1-a}} \|\vp_\circ\|_r\,[t^\frac
1p-(t-s)^\frac 1p]^{-\frac{a}{1-a}},$$ whence
\be\label{estb2aa}\|\vp(s)\|_2\leq c M^\frac{(2-p)\gamma}{2}\|\vp_\circ\|_r\, [t^\frac 1p-(t-s)^\frac
1p]^{-p\gamma },\ee which, by the expression of $M$,
 gives \eqref{Maxpostaa}.
\par Assume that $b=0$. Then, from Lemma\,\ref{adgt},
estimate \eqref{max11a} holds for all $r\in [1,2]$.
The proof of \eqref{Maxpost} follows step by step the above proof.
 One has just to use \eqref{maxa3aa} with $r=1$:
  $$\|\vp\|_2\leq c\|\nabla
\vp\|_p^a\,\|\vp\|_1^{1-a}\leq c\|\nabla
\vp\|_p^a\,\|\vp_\circ\|_1^{1-a}\,,$$ where now $$
a=\mbox{\large$\frac{np}{2(np+p-n)}$}\,,\qquad 1-a=\mbox{\large$\frac{np+2p-2n}{2(np+p-n)}$}\,.
$$ Note that $a<1$ if and only if $p>\frac{2n}{n+2}$. With the previous calculations we arrive at
\eqref{estb2aa} with $r=1$, which, by setting
$\gamma=\gamma(1)=\beta$ gives \be\label{estb2}\|\vp(s)\|_2\leq c\, M^\frac{(2-p)\beta}{2}\|\vp_\circ\|_1\,[t^\frac 1p-(t-s)^\frac
1p]^{-p\beta},\ee hence \eqref{Maxpost} for $r=2$.
 Finally, using the $L^p$-convexity inequality, then \eqref{estb2} and \eqref{max11a} with $r=1$, we get
 $$ \|\vp(s)\|_r\leq \|\vp(s)\|_2^{\frac{2}{r'}}\|\vp(s)\|_1^{\frac{r'-2}{r'}}\leq  c\,M^\frac{(2-p)\beta}{r'}\|\vp_\circ\|_1 [t^\frac 1p-(t-s)^\frac
1p]^{-\frac{2p\beta}{r'}}$$
 which gives estimate \eqref{Maxpost},
  and completes the proof.\chiu

 \begin{lemma}\label{compattificazione}{\sl
 Let $\nabla\psi\in L^2((0,t)\times \O)$, $\nabla v \in L^r((0,t)
 \times\O)$, for some $r>1$, and let $h^m$ be a
 sequence with $\nabla h^{m}$ bounded in $ L^2((0,t)
 \times \O)$, uniformly in $m\in \N$. Then, there exists a subsequence $h^{m_k}$ such
 that
 $$\lim_{k\to \infty}\int_0^s \int_\O\left((\mu+
 |J_{\frac {1}{m_k}}(\nabla v)|^2)^\frac{p-2}{2}-
 (\mu+|\nabla v|^2)^\frac{p-2}{2}\right)
\nabla h^{m_k}\cdot  \nabla
\psi \,dx\,d\tau=0\,.$$ }
 \end{lemma}
 \Pr
 By assumption  $\nabla v\in L^{p}((0,t)\times \O)$. Therefore
one has $J_{\frac 1m}(\nabla v) \to \nabla v$ in $L^p((0,t)\times
\O)$ as $m$ goes to $\infty$. This ensure the existence of a
subsequence $ J_{\frac {1}{m_k}}(\nabla v)$ converging  to $\nabla
v$ a.e. in $(t,x)\in (0,T)\times \O$. Therefore, along this
subsequence,
$$(\mu+|J_{\frac {1}{m_k}}(\nabla
v)|^2)^\frac{p-2}{2}\to (\mu+|\nabla v|^2)^\frac{p-2}{2},\ \textrm{
a.e. in } (s, x)\in (0,t)\times \O\,.$$
 By H\"{o}lder's inequality  one has
 $$\ba{ll}\dy \int_0^s\!\! \int_\O\left((\mu+|J_{\frac {1}{m_k}}
  (\nabla v)|^2)^\frac{p-2}{2}-(\mu+|\nabla v|^2)^\frac{p-2}{2}\right)
\nabla h^{m_k}\cdot  \nabla \psi \,dx\,d\tau\\ \dy \leq
\left(\int_0^s\!\!\!\int_\O\! |\nabla \psi
(\tau)|^2\!\left((\mu+|J_{\frac {1}{m_k}}(\nabla
v)|^2)^\frac{p-2}{2}-(\mu+|\nabla
v|^2)^\frac{p-2}{2}\right)^{\!2}dxd\tau\!\right)^{\!\!\frac
12}\!\|\nabla h^{m_k}\|_{L^2((0,t)\times \O)}.\ea$$ Further, since
$\nabla \psi\in L^{2}((0,t)\times (\O))$, we have
$$\ba{ll}\dy\vs1 |\nabla
\psi(\tau)|^2\left((\mu+|J_{\frac {1}{m_k}}(\nabla v)|^2)^\frac{p-2}{2}-(\mu+|\nabla v|^2)^\frac{p-2}{2}\right)^2
\\ \hskip4cm \leq
\left(2\mu^\frac{p-2}{2}\right)^2|\nabla\psi(\tau)|^2 \in
L^{1}(0,t;L^1(\O))\,.\ea$$ Therefore, from the Lebesgue dominated
convergence theorem we obtain the result. \chiu

  \begin{lemma}\label{semigr}{\sl Let $\nu>0$,
   $\mu>0$ and  $p >\frac{2n}{n+2}$. Let $v(x)$ satisfy
\be\label{beeta2st}
\,\|(\mu+|\nabla v|^2)^\frac 12\|^p_{p} \leq
 M\,,\ee
with a positive constant $M$. Then, for any $\vp_\circ(x)\in
C_0^{\infty}(\O)$, there exists a unique solution $\vp\in C(0,t;
L^{2}(\O))\cap L^2(0,t; W_0^{1,2}(\O))$ of
 the following integral equation
\be\label{AD1wf}\begin{array}{l}\vs1 \dy
\int_0^s(\vp,\psi_\tau)d\tau-\nu\int_0^s (\nabla\vp,
\nabla\psi)\,d\tau -\int_0^s \left((\mu+|\nabla
v|^2)^\frac{p-2}{2}\nabla \vp, \nabla
\psi \right)\,d\tau\\
\hskip4cm = (\vp(s),
\psi(s))-(\vp_\circ,\psi(0))\,,\
\forall
\psi\!\in\!C_0^\infty([0,t)\times
\O) .\end{array}\ee Moreover, for
any $r\in [1,2]$, one has
\be\label{server}\|\vp(s)\|_r\leq
\|\vp_\circ\|_r,\ \forall s \in
[0,t]\,,\ee
\be\label{Maxpostnonu}\|\vp(s)\|_r\leq
c\,M^\frac{(2-p)\beta}{r'}\|\vp_\circ\|_1\,
s^{-{p\beta}\left(1-\frac{1}{r}\right)}\,,\quad
\forall s\in(0,t]\,,\ee and
\be\label{maxgrad} \|\nabla
\vp(s)\|_p\leq
c\,M^{\frac{2-p}{2}(\beta+\frac
1p)}\|\vp_\circ\|_1\,
s^{-\frac{\beta p+1}{2}} \,,\quad
\forall s\in(0,t]\,, \ee
 with $\beta$ given in \eqref{Maxpostr}  and $r'$ conjugate exponent of $r$.
  }\end{lemma}
\Pr Let us consider the unique solution of system \eqref{AD1} given in Lemma\,\ref{newag1}.
 The solution $\vp=\vp(\eta)$ satisfies estimates \eqref{server} and \eqref{Maxpostnonu},
  uniformly in $\eta>0$. The proof is the same of Lemma\,\ref{newag2}in the case where $b=0$, replacing the assumption \eqref{beeta2} by \eqref{beeta2st}.
Let us multiply \eqref{AD1} by $\vp_s$ and integrate over $\O$. We
get
 \be\label{grad1d}\frac \nu2\frac{d}{ds}
\|\nabla\vp\|_2^2+\frac 12\frac{d}{ds}
 \|(\mu+|J_\eta(\nabla v)|^2)^\frac{p-2}{4}\nabla\vp\|_2^2+\|\vp_s\|_2^2=0\,.\ee
 Let us multiply \eqref{AD1} by $\vp$ and integrate over $\O$.
By using estimate \eqref{Maxpostnonu}  and recalling \eqref{amu} we
get
 $$\nu\|\nabla\vp\|_2^2+ \|a_\eta(\mu,v)^\frac 12\nabla\vp\|_2^2\leq \|\vp_s\|_2\|\vp\|_2
 \leq   c\|\vp_s\|_2\|\vp_\circ\|_1\, M^\frac{(2-p)\beta}{2}s^{-\frac{\beta p}{2}}\,,\  \forall s\in(0,t]\,,
 $$
 hence
  $$c\frac{(\,\nu\|\nabla \vp\|_2^2+  \|a_\eta(\mu,v)^\frac 12\nabla\vp\|_2^2\,)^2}{
  \|\vp_\circ\|_1^2\, M^{(2-p)\beta}}\, s^{{\beta p}}\leq  \|\vp_s\|_2^2\,,\quad \forall s\in(0,t]\,.
 $$
 Using this estimate in \eqref{grad1d} we get
 $$\frac 12\frac{d}{ds}\left(\,\nu\|\nabla\vp\|_2^2+
 \|a_\eta(\mu, v)^\frac{1}{2}\nabla\vp\|_2^2\,\right)+c
  \frac{ (\nu\|\nabla\vp\|_2^2+ \|a_\eta(\mu, v)^\frac{1}{2}\nabla\vp\|_2^2)^2}{\|\vp_\circ\|_1^2\,
  M^{(2-p)\beta}}\, s^{{\beta p}}\leq 0\,,
  $$
  which, integrated from $0$ to $s$ gives
  $$ \nu\|\nabla \vp\|_2^2+\|a_\eta(\mu,v)^\frac 12\nabla\vp\|_2^2\leq c\, M^{(2-p)\beta}  \|\vp_\circ\|_1^2\, s^{{-(\beta p+1)}}\,. $$
  By using the arguments in \eqref{maxa1}, the above estimate and assumption \eqref{beeta2} on $v$,  we have
$$\ba{ll}\dy\vs1 \|\nabla \vp\|_p\!\!\!\!&\leq \dy
\|(\mu+|J_\eta(\nabla v)|^2)^\frac{(p-2)}{4}\,\nabla \vp\|_2\,
\|(\mu+|J_\eta(\nabla
v)|^2)\|_{\frac p2}^\frac{(2-p)}{4}\\
\dy &\leq\dy c\,
M^{\frac{(2-p)\beta}{2}}M^{\frac{2-p}{2p}}
\|\vp_\circ\|_1
\, s^{{-\frac{\beta
p+1}{2}}}\,,\ \forall s\in(0,t],\ea$$ which gives \eqref{maxgrad}.
In order to obtain the result, the next step is to prove that,
denoting for any $\eta>0$ by $\vp^\eta$ the unique solution of
\eqref{AD1}, which satisfies estimates \eqref{server},
\eqref{Maxpostnonu} and \eqref{maxgrad}, the sequence $\{\vp^\eta\}$
converges in some sense to the solution of the integral equation
\eqref{AD1wf} as $\eta$ goes to zero. The proof is straightforward.
We avoid the details and just show the following  convergence, along
a suitable subsequence,
$$\int_0^s \left((\mu+|J_\eta(\nabla v)|^2)^\frac{p-2}{2}\nabla \vp^\eta, \nabla
\psi \right)\,d\tau\to \int_0^s \left((\mu+|\nabla
v|^2)^\frac{p-2}{2}\nabla \vp, \nabla \psi \right)\,d\tau\,.  $$
Indeed, writing
$$\ba{ll}\dy\vs1 \int_0^s \left((\mu+|J_\eta(\nabla v)|^2)^\frac{p-2}{2}\nabla \vp^\eta, \nabla
\psi \right)\,d\tau-\int_0^s \left((\mu+|\nabla
v|^2)^\frac{p-2}{2}\nabla \vp, \nabla \psi \right)\,d\tau\,\\
\vs1\dy =\int_0^s\int_\O\left((\mu+|J_\eta(\nabla
v)|^2)^\frac{p-2}{2}-(\mu+|\nabla v|^2)^\frac{p-2}{2}\right) \nabla
\vp^\eta\cdot  \nabla \psi \,dx\,d\tau\\ \dy -\int_0^s\int_\O
(\mu+|\nabla v|^2)^\frac{p-2}{2} (\nabla \vp^\eta-\nabla \vp)\cdot
\nabla \psi\,dx \,d\tau, \ea
$$
the first integral goes to zero, thanks to Lemma\,\ref{compattificazione}
and the second integral goes to zero thanks to the weak convergence of $\nabla\vp^\eta$ to $\nabla\vp$
 in $L^2((0,t)\times \O)$ (using \eqref{max5}).    \chiu
 \vskip0.2cm
We also give a useful inequality, referring, for instance, to \cite{CG}.
\begin{lemma}\label{lCarlo2}{\sl
Let  $\mu>0$. For any given real numbers
$\xi,\eta\geq 0$, the following inequality holds:
$$\bigg|\frac1{(\mu+\xi)}{\!\atop
^{{2-p}}}-\frac1{(\mu+\eta)}{\!\atop ^{{2-p}}}
\bigg|\le\frac{2-p}{\mu}{\!\atop
^{{\!\!\!\!\!\!3-p}}}\,|\xi-\eta|.$$}
\end{lemma}
\vskip0.2cm

Below we recall some well known results for Bochner spaces.
\begin{lemma}\label{solonnl}{\sl
Let $u$ belong to $L^\infty(\ve,T; W^{2,\ov q}(\O)\cap W_0^{1,\ov q}(\O))$
 with $u_t\in L^\infty(\ve, T; L^{\ov q}(\O))$. For $m=0,1$, if $\lambda =2-m-\frac {n}{\ov q}\in (0,1)$, then
$$[\,\nabla^m u \,]_{\lambda,t,x}\leq C \big[\sup_{(\ve,T)} (\|u_t(t)\|_{\ov q}+ \|D^2 u(t)\|_{\ov q})
+\sup_{(\ve,T)}\|u(t)\|_{\ov q}\,\big]\,, $$with $C$ independent of $u$.
}\end{lemma}
\Pr The proof is a trivial generalization of Theorem 2.1 proved in \cite{Sol77}. \chiu\vskip0.2cm
For the following embedding results we refer, for instance, to \cite{show}, Ch. 3.
\begin{lemma}\label{boch}{\sl Let $X$ be a Banach space and let $X'$ be its dual.
Assume that $X$ is dense and
continuously embedded in an Hilbert space $H$. We identify $H$ with $H'$,
which is continuously embedded in $X'$. If $u\in L^q(0,T;X)$ and $u'\in L^{q'}(0,T;X')$,
 with $q, q'\in (1,+\infty)$, $\frac 1q+\frac{1}{q'}=1$,
 then $u$ is almost everywhere equal to a continuous function from $[0,T]$ into $H$.  }
\end{lemma}

\begin{lemma}\label{aubin}{\rm [Aubin-Lions]}- {\sl Let $X$, $X_1$, $X_2$ be
Banach spaces. Assume that   $X_1$ is compactly embedded in $X$ and $X$
 is continuously embedded in $X_2$, and that $X_1$ and $X_2$ are reflexive.
 For $1<q,s<\infty$, set $$W=\{\psi \in L^s(0,T; X_1): \psi_t\in L^q(0,T; X_2)\}\,.$$
Then the inclusion $W \subset L^s(0,T; X)$ is compact.} \end{lemma}

\section*{\normalsize 3. Approximating systems}
\renewcommand{\theequation}{3.\arabic{equation}}
\renewcommand{\thetho}{3.\arabic{tho}}
\renewcommand{\thedefi}{3.\arabic{defi}}
\renewcommand{\therem}{3.\arabic{rem}}
\renewcommand{\theprop}{3.\arabic{prop}}
\renewcommand{\thelemma}{3.\arabic{lemma}}
\renewcommand{\thecoro}{3.\arabic{coro}}
\setcounter{equation}{0} \setcounter{coro}{0} \setcounter{lemma}{0}
\setcounter{defi}{0} \setcounter{prop}{0} \setcounter{rem}{0}
\setcounter{tho}{0}
 Let us study the approximating systems
\be\label{PFl}\begin{array}{ll}\dy
u_t-\nabla\cdot\left(\left(\mu+|\nabla
u|^2\right)^\frac{(p-2)}{2}\nabla u\right)= 0\,,&\hskip-0.2cm\textrm{ in
}(0,T)\times\O,
\\\dy \vs1\hskip4cmu(t,x)=0\,,&\hskip-0.2cm\textrm{ on }
(0,T)\times\po,\;\\\dy\hskip4cm u(0,x)=u_\circ(x),&\hskip-0.2cm\mbox{ on
}\{0\}\times\O,\end{array}\ee with $\mu>0$, and

\be\label{PFepv}\begin{array}{ll}\dy v_t-\nu\Delta
v-\nabla\cdot\left(\left(\mu+|\nabla
v|^2\right)^\frac{(p-2)}{2}\nabla v\right)= 0\,,&\hskip-0.2cm\textrm{ in
}(0,T)\times\O,
\\\dy\vs1\hskip5cm v(t,x)=0\,,&\hskip-0.2cm \textrm{ on }
(0,T)\times\po,\;\\\dy\hskip5cm v(0,x)=v_\circ(x),&\hskip-0.2cm\mbox{ on
}\{0\}\times\O\,,\end{array}\ee %
with $\mu>0$ and $\nu>0$.\par Let us introduce the operators from $V$ to $V'$ defined as
\be\label{oper1}A(v):=-\nabla\cdot\left(|\nabla v|^{p-2}\nabla
v\right)\,,\ee \be\label{oper2}A_\mu(v):=-
\nabla\cdot\left(\left(\mu+|\nabla v|^2\right)^\frac{(p-2)}{2}\nabla
v\right)\,.\ee
They are both
monotonous and emicontinuous
operators\footnote{For the sake of
brevity for the corresponding
definitions we refer, for instance,
to \cite{lions} Ch. II, Sec.1.2.}.
 Set \be\label{defB}B(\mu, w):= c\,\|w\|_2^2+
c(\O,T)\mu^\frac p2\,. \ee Recall that, from \eqref{oper3},
$$a(\mu,v)=\left(\mu+|\nabla v|^2\right)^\frac{(p-2)}{2}\,.$$
\begin{defi}\label{wmu} {\rm Let $\mu> 0$. Let $u_{\circ}\in
L^2(\O)$. A field $u\!:(0,T)\times \O\to\R^N$
 is said a solution of
system {\rm \eqref{PFl}} if
$$
\ba{ll}\vs1\dy u\in  L^{p}(0,T; V)\cap C(0,T;L^2(\O))\,, \ t^\frac 1p \nabla u\in L^{\infty}(0,T;L^p(\O))\,,\\
\dy \ t\,u_t\in L^{\infty}(0,T;L^2(\O))\,,\ t^\frac{p+2}{2p}\,\nabla
u_t \in L^{2}(0,T; L^p(\O))\,,\ea$$
$$\ba{ll}\dy \vs1\int_0^t \left[(u,\psi_\tau) -\left(a(\mu,
u)\,\nabla u,\nabla \psi\right)\right]d\tau=(u(t),
\psi(t))-(u_\circ, \psi(0)),\\ \hskip5cm   \forall \psi\in
W^{1,2}(0,T;L^2(\O))\cap L^p(0,T;V),\ea$$ and $$\lim_{t\to
0^+}\|u(t)-u_\circ\|_2=0\,.$$}
\end{defi}
\begin{defi}\label{wmunu}
{\rm Let $\mu> 0$, $\nu>0$. Let $v_{\circ}\in L^2(\O)$. A field
$v\!:(0,T)\times\O\to\R^N$
 is said a  solution of
system {\rm \eqref{PFepv}} if $$ v\in C(0,T;L^2(\O))\cap L^{2}(0,T;
W_0^{1,2}(\O)),$$ $$ v_t\in L^{\infty}(0,T;L^2(\O))\cap
L^2(0,T;W^{1,2}(\O)),$$ $$\ba{ll}\dy\vs1 \int_0^t
\left[(v,\psi_\tau)-\nu(\nabla v,\nabla\psi)-\left(a(\mu, v)\,\nabla
v,\nabla \psi\right)\right]d\tau=(v(t),\psi(t))-(v_\circ,
\psi(0)),\\ \hfill \forall \psi \in W^{1,2}(0,T;L^2(\O))\cap
L^2(0,T;W_0^{1,2}(\O))\,,\ea$$ and $$\lim_{t\to
0^+}\|v(t)-v_\circ\|_2=0\,.$$}
\end{defi}

We have introduced the
approximating systems \eqref{PFl}
and \eqref{PFepv} in order to prove
Theorem\,\ref{existencenomu}. In
particular, the introduction of
this second kind of approximation
is necessary, for our technique, to
obtain the crucial estimate $u_t\in
L^\infty(\ve,T; L^q(\O))$, which is
one of the key tools to get higher
order integrability for the second
derivatives (Theorem\,\ref{regularitynomu} and Theorem\,\ref{regularitynomun2}) and,
further, space-time H\"{o}lder
regularity (see Corollary\,\ref{corregularitynomu}).
 If we limit ourselves just
 to Theorem\,\ref{existencenomu},
 where the $L^2(\ve,T; W^{2,p}(\O))$
  integrability of the second
  derivatives is shown, then we could
 avoid the study of system \eqref{PFepv},
 making the proof easier.\par
In order to study the existence and regularity of a solution of
 \eqref{PFl}, firstly we study the same issues for
the parabolic approximating system \eqref{PFepv}. The existence
and regularities for the solution of this latter system are obtained
 in the following Proposition\,\ref{existence},
  Corollary\,\ref{existencecor1} and Corollary\,\ref{pregcor}. These results are a fundamental
 step for the proof of our main results. On the other hand their
 proofs rely on the Galerkin approximation method, with a suitable
 choice of the basis functions, and related weighted estimates
 in $W^{m,r}$-spaces, $m \in \N\cup\{0\}$.  Since these arguments are
 standard, we confine the proofs in the appendix.
We observe that in the propositions
 below we will assume that the initial
 data of problems \eqref{PFl} and
 \eqref{PFepv} are in $C_0^{\infty}(\O)$.
 This assumption could be weakened,
 for the validity of the same results.
 This will be done in the next sections,
 where we deal with the solutions of problem \eqref{PF}
  and of problem \eqref{PFl}, and consider the
  completion of
  $C_0^\infty(\O)$ in  $L^2(\O)$ and
  $W_0^{1,2}(\O)$,
  so realizing a suitable generalization of the
  following results.

\begin{prop}\label{existence}
{\sl Let be $\nu>0$ and $\mu>0$. Assume that $v_\circ$
belongs to $C_0^\infty(\O)$. Then there exists a unique solution $v$
of system \eqref{PFepv} in the sense of Definition\,\ref{wmunu}.
In particular:\begin{itemize}
\item[i)]
 $v\in L^{\infty}(0,T;L^2(\O))\cap  L^{p}(0,T; V)\,$, uniformly in $\nu$
and $\mu$; \item[ii)]
$ v_t\in
L^{\infty}(0,T;L^2(\O))$, $\,\nabla v_t \in L^{2}(0,T; L^p(\O))$, uniformly in $\nu$, non-uniformly in $\mu$;
\item[iii)] For all $\nu>0, v\!\in\!L^2(0,T; W_0^{1,2}(\O)),
v_t \in L^{2}(0,T; W^{-1,2}(\O)), \nabla v_t\in L^2(0,T;
L^2(\O))$.
\end{itemize}}
\end{prop}
\begin{coro}
\label{existencecor1} {\sl Under the assumptions of Proposition\,\ref{existence}, we have
\begin{itemize}
\item[i)]
$t^\frac 1p \nabla v\in L^{\infty}(0,T;L^p(\O))$,
 $t\,v_t\in
L^{\infty}(0,T;L^2(\O))$,
$t^\frac 12\, v_t\!\in L^2(0,T;L^2(\O))$, $t^\frac{p+2}{2p}\,\nabla v_t \!\in
L^{2}(0,T; L^p(\O))$,  uniformly in $\nu>0$ and $\mu>0$;
\item[ii)]  $t\,\nabla
v_t\in L^2(0,T; L^2(\O))$, non-uniformly in $\nu>0$\,.
\end{itemize}}
\end{coro}
 Set
 \be\label{alpha} \alpha:=\mbox{\large$ \frac{4}{2p-n(2-p)}\,,$}\ee
 \be\label{beta12}
 \beta_1(p):=\mbox{\large$\frac12$}\left[\alpha -\mbox{\large$\frac{n(2-p)^2}{p(2p-n(2-p))}$}\right],\ \
\beta_2(p):=\mbox{\large$\frac 12$}\left[\alpha+\mbox{\large$\frac{2-p}{p}$}\right]\,.\ee\begin{coro}\label{pregcor} { Let $\mu> 0$. Assume that $v_\circ$ belongs to
$C_0^\infty(\O)$. If $p>\frac 32$, then  $\nabla v\in
C(0,T;L^{2}(\O))$, $v\in
L^{2}(0,T;W^{2,p}(\O))$, with
$$ \|\nabla v\|_{C(0,T; L^2(\O))}+\|v\|_{L^2(0,T; W^{2,p}(\O))}\leq M_3(\|v_\circ\|_{1,2}, B(\mu, v_\circ))\,.$$
 If $p>p_\circ$, then $t^{\beta_1(p)}\, \nabla v\in
C(0,T;L^{2}(\O))$, $t^{\beta_2(p)}\, v\in
L^{2}(0,T;W^{2,p}(\O))$,  with
$$\ba{ll}\vs1\dy t^\alpha \|\nabla
v(t)\|_2^2\leq C\, \mu^\frac{n(2-p)^2}{2(2p-n(2-p))}\,B(\mu,
v_\circ)^\frac 2p\, t^{\frac{n(2-p)^2}{p(2p-n(2-p))}}\\
\hskip2cm\dy+C\,(B(\mu, v_\circ))^\frac{4-n(2-p)}{2p-n(2-p)} +c
B(\mu, v_\circ) t^\alpha(1+t)\,,\ea$$
$$\ba{ll}\dy\int_0^t\tau^{\alpha+\frac{2-p}{p}} \|D^2
v(\tau)\|_p^2\, d\tau
 \leq
C(B(\mu, v_\circ),T)\,,\ea$$
where $\alpha$, $\beta_1$ and $\beta_2$ are given by \eqref{alpha}--\eqref{beta12}. }
\end{coro}
In the next proposition, starting from the
 existence and regularities of the solution of
  system \eqref{PFepv}, given in Proposition\,\ref{existence},
 and passing to the limit as $\nu$ goes to zero,
 we deduce analogous existence and regularity
properties for the solution of system \eqref{PFl}.
\begin{prop}\label{existencenonu}
{\sl Let be $\mu>0$. Assume that $u_\circ$ belongs
to $C_0^\infty(\O)$. Then there exists a unique solution $u$ of
system \eqref{PFl} in the sense of Definition\,\ref{wmu}, such that
\begin{itemize}
\item[i)] $u\in L^\infty(0,T; L^2(\O))\cap L^2(0,T; V)$, uniformly in $\mu$;
\item[ii)] for all $\mu>0$, $
u_t\in  L^{\infty}(0,T;L^2(\O))$, $\,\nabla u_t \in L^{2}(0,T;
L^p(\O))$.
\end{itemize}
}
\end{prop}
\Pr From Proposition\,\ref{existence} for all $\mu>0$, the function $v$,
solution of \eqref{PFepv} corresponding to the initial data
$u_\circ\in C_0^\infty(\O)$,  satisfies   the
bounds collected below, uniformly in $\nu>0$,
$$\ba{ll}\dy\vs1
\|v^\nu \|_{L^{\infty}(0,T;L^2(\O))}+\|v^\nu
\|_{L^{p}(0,T;V)}\leq M(T,\O,
\|u_\circ\|_{2})\,,\ea$$
$$ \|v^\nu_t\|_{L^{\infty}(0,T;L^2(\O))}+\|v^\nu _t\|_{L^{2}(0,T;
W^{1,p}(\O))}+\|v^\nu
_t\|_{L^{p'}(0,T;V')}\leq
M_1(\mu,T,\O,
\|u_\circ\|_{2,2})\,,$$ where the
constant $M_1$ blows up as $\mu\to
0$. Hence we can extract a
subsequence, still denoted by
$\{v^\nu\}$, weakly or weakly-*
converging,  in the above norms, to
a function $u$, as $\nu$ tends to
zero. Further, from $v^\nu\in
L^p(0,T; V)$, it follows that
$a(\mu, v^\nu)\,\nabla v^\nu \in
L^{p'}(0,T; L^{p'}(\O))$, hence
$A_\mu(v^\nu) \in L^{p'}(0,T; V')$
and (see \eqref{numA})
\be\label{b5}A_\mu(v^\nu)\rightharpoonup
\widetilde\chi\ \textrm{in }
L^{p'}(0,T; V')\
\textrm{weakly}.\ee Using the
monotonicity trick as in
\cite{lions}, we show that the
non-linear part $A_\mu(v^\nu)$
actually (weakly) converges to
$A_\mu(u)$ and that the limit $u$
is a solution of system
\eqref{PFl}. In particular the
regularities stated for $u$ follow
from the analogous regularities of
$v^\nu$ and the lower
semi-continuity of the norm for the
weak convergence.

 Hence, set
$$X_{\nu}=\int_0^t(A_\mu(v^\nu)-A_\mu(\vp), v^\nu-\vp)d\tau \geq 0,\ \forall
\vp\,\in L^p(0,T;V)\,.$$ By using that $v^\nu$ is a solution of
\eqref{PFepv},  we can write $X_\nu$ as follows
$$\ba{ll}\vs1\dy X_{\nu}\!\!&=\dy\int_0^t(A_\mu(v^\nu),v^\nu)d\tau-
\int_0^t(A_\mu(v^\nu),\vp)d\tau -\int_0^t(A_\mu(\vp),
v^\nu-\vp)d\tau\\\vs1\dy &=\dy\frac 12 \|v^\nu(0)\|_2^2-\frac
12\|v^\nu(t)\|_2^2-\nu\int_0^t\|\nabla v^\nu(\tau)\|_2^2d\tau\\&\dy\quad -
\int_0^t(A_\mu(v^\nu),\vp)d\tau -\int_0^t(A_\mu(\vp),
v^\nu-\vp)d\tau\,.\dy \ea
$$
Passing to the lim sup and observing that $$\lim
\sup(-\nu\int_0^t\|\nabla v^{\nu}(\tau)\|_2^2d\tau)=-\lim\inf
(\nu\int_0^t\|\nabla v^{\nu}(\tau)\|_2^2)d\tau\leq 0,$$  we get
$$\ba{ll}\dy 0\leq \lim \sup X_{\nu}\leq \frac 12\|u_\circ\|_2-\frac 12\|u(t)\|_2^2
+\lim \sup(-\nu\int_0^t\|\nabla v^{\nu}(\tau)\|_2^2d\tau)
\\\hskip3cm
\dy - \int_0^t(\widetilde\chi,\vp)d\tau -\int_0^t(A_\mu(\vp),
u-\vp)d\tau\\\dy \quad\leq \frac 12\|u_\circ\|_2^2-\frac 12\|u(t)\|_2^2-
\int_0^t(\widetilde\chi,\vp)d\tau -\int_0^t(A_\mu(\vp), u-\vp)d\tau.
\ea
$$
On the other hand, from the weak convergence it is easy to see that
the limit $u$ satisfies \be\label{appp2} \int_0^t(\widetilde\chi, u)
d\tau =\frac 12\|u_\circ\|_2^2-\frac 12\|u(t)\|_2^2\,.\ee
 Therefore

$$\int_0^t(\widetilde\chi-A_\mu(\vp), u-\vp)d\tau \geq 0. $$
Taking $\vp=u+\lambda w$, for $\lambda>0$ and for some $w\in
L^{p}(0,T, V)$, and then letting $\lambda$ tend to zero the thesis
follows.\chiu \vskip0.5cm%
Next, we improve property {\it ii)} of Proposition\,\ref{existencenonu}, by using Corollary\,\ref{existencecor1} and Corollary\,\ref{pregcor}. Indeed we are going to show that time weighted estimates for $u_t$ and  $\nabla u_t$ hold for any $\mu>0$.
\begin{prop}\label{existpeso1}{\sl
Let be $\mu>0$. Assume that $u_\circ$ belongs
to $C_0^\infty(\O)$. Then the solution $u$ of system \eqref{PFl} of Proposition\,\ref{existencenonu} satisfies,
uniformly in $\mu$, $\,t^\frac 1p \nabla u\in
L^{\infty}(0,T;L^p(\O))$
 $t\,u_t\in
L^{\infty}(0,T;L^2(\O))$, $t^\frac{p+2}{2p}\,\nabla u_t \in
L^{2}(0,T; L^p(\O))$. \par Further, for $p>\frac 32$, $\nabla
u\in C(0,T;L^{2}(\O))$, $u\in L^{2}(0,T;W^{2,p}(\O))$, with$$
\|\nabla u\|_{C(0,T; L^2(\O))}+ \|u\|_{L^2(0,T; W^{2,p}(\O))}\leq
M_3(\|u_\circ\|_{1,2}, B(\mu, u_\circ))\,.$$ For $p>p_\circ$,
$t^{\beta_1(p)}\, \nabla u\in C(0,T;L^{2}(\O))$, $t^{\beta_2(p)}\,
u\in L^{2}(0,T;W^{2,p}(\O))$,  with
$$\ba{ll}\vs1\dy t^\alpha \|\nabla
u(t)\|_2^2\leq C\, \mu^\frac{n(2-p)^2}{2(2p-n(2-p))}\,B^\frac
2p(\mu,
u_\circ)\, t^{\frac{n(2-p)^2}{p(2p-n(2-p))}}\\
\dy \hskip2cm +C\,(B(\mu, u_\circ))^\frac{4-n(2-p)}{2p-n(2-p)} +c
B(\mu, u_\circ) t^\alpha(1+t)\,,\ea$$
$$\ba{ll}\dy\int_0^t\tau^{\alpha+\frac{2-p}{p}} \|D^2
u(\tau)\|_p^2\, d\tau
 \leq
C(B(\mu, u_\circ),T),\ea$$ where $\alpha$, $\beta_1$ and $\beta_2$ are given by \eqref{alpha}--\eqref{beta12}. }
\end{prop}
\Pr Let us consider system \eqref{PFepv}, with initial data
$v_\circ=u_\circ$. From Corollary\,\ref{existencecor1},
$$ \|t^\frac 1p\,\nabla
v^\nu\|_{L^{\infty}(0,T;L^p(\O))}+\|t\,v^\nu_t\|_{L^{\infty}(0,T;L^2(\O))}+\|t^\frac{p+2}{2p}\nabla
v^\nu_t\|_{L^{2}(0,T; L^p(\O))} \leq M_2(T,\O, \|u_\circ\|_{2})$$
and, from Corollary\,\ref{pregcor}, for $p>\frac 32$
$$\|\, \nabla v^\nu\|_{C(0,T;L^{2}(\O))}+
\|\, v^\nu\|_{L^{2}(0,T;W^{2,p}(\O))}\leq M_3( \|u_\circ\|_{1,2},
B(\mu, u_\circ))\,,$$ while for $p>p_\circ$
$$\ba{ll}\vs1\dy t^\alpha \|\nabla
v^\nu(t)\|_2^2\leq C\, \mu^\frac{n(2-p)^2}{2(2p-n(2-p))}\,B(\mu,
u_\circ)^\frac 2p\, t^{\frac{n(2-p)^2}{p(2p-n(2-p))}}\\
\dy \hskip2cm +C\,(B(\mu, u_\circ))^\frac{4-n(2-p)}{2p-n(2-p)} +c
B(\mu, u_\circ) t^\alpha(1+t)\,,\ea$$
$$\int_0^t\tau^{\alpha+\frac{2-p}{p}} \|D^2 v^\nu(\tau)\|_p^2\,
d\tau
 \leq
C(B(\mu, u_\circ),T),$$ uniformly in $\nu>0$. Therefore, passing to
the limit as $\nu$ tends to zero, and then reasoning as in the proof
of Proposition\,\ref{existencenonu} we get that the limit $u$,
solution of system
 \eqref{PFl}, satisfies the same bounds. The proof is
then completed.
\chiu

\section*{\normalsize 4. Proof of Theorem\,\ref{existencenomu}
 and Corollary\,\ref{corexistencenomu}}
\renewcommand{\theequation}{4.\arabic{equation}}
\renewcommand{\thetho}{4.\arabic{tho}}
\renewcommand{\thedefi}{4.\arabic{defi}}
\renewcommand{\therem}{4.\arabic{rem}}
\renewcommand{\theprop}{4.\arabic{prop}}
\renewcommand{\thelemma}{4.\arabic{lemma}}
\renewcommand{\thecoro}{4.\arabic{coro}}
\setcounter{equation}{0} \setcounter{coro}{0} \setcounter{lemma}{0}
\setcounter{defi}{0} \setcounter{prop}{0} \setcounter{rem}{0}
\setcounter{tho}{0}
In the next proposition, under the assumption of an initial data $u_\circ\in L^2(\O)$ or $u_\circ\in W_0^{1,2}(\O)$,
  we study the existence and regularity of the solution
 of system \eqref{PFl}, for any $\mu>0$, already obtained in Proposition\,\ref{existencenonu}
  and Proposition\,\ref{existpeso1} under the stronger assumption of $u_\circ\in C_0^\infty(\O)$.
   \begin{prop}\label{DDAAGG}
{\sl Let be  $\mu>0$. Assume that $u_\circ$ belongs
to $L^2(\O)$. Then there exists a unique solution $u$ of
system \eqref{PFl} in the sense of Definition\,\ref{wmu}.
Moreover, for $p>p_\circ$, $t^{\beta_1(p)}\, \nabla u\in
C(0,T;L^{2}(\O))$, $t^{\beta_2(p)}\, u\in
L^{2}(0,T;W^{2,p}(\O))$,  with
$$\ba{ll}\vs1\dy t^\alpha \|\nabla
u(t)\|_2^2\leq C\, \mu^\frac{n(2-p)^2}{2(2p-n(2-p))}\,B^\frac
2p(\mu, u_\circ)\, t^{\frac{n(2-p)^2}{p(2p-n(2-p))}}\\
\dy \hskip2cm +C\,(B(\mu, u_\circ))^\frac{4-n(2-p)}{2p-n(2-p)} +c
B(\mu, u_\circ) t^\alpha(1+t)\,,\ea$$
$$\ba{ll}\dy\int_0^t\tau^{\alpha+\frac{2-p}{p}} \|D^2
u(\tau)\|_p^2\, d\tau
 \leq
C(B(\mu, u_\circ),T),\ea$$
where $\alpha$, $\beta_1$ and $\beta_2$ are given by \eqref{alpha}--\eqref{beta12}.
Finally, assume that $u_\circ\in
W_0^{1,2}(\O)$. Then, for $p>\frac 32$, $\nabla u\in
C(0,T;L^{2}(\O))$, $u\in L^{2}(0,T;W^{2,p}(\O))$, with
$$ \|\nabla u\|_{C(0,T; L^2(\O))}+\|u\|_{L^2(0,T; W^{2,p}(\O))}\leq M_3(\|u_\circ\|_{1,2}, B(\mu,
u_\circ))\,.$$ }\end{prop}

\Pr
 Let $\dy\{ u_\circ^m(x)\}$ be a sequence
in $C_0^\infty(\O)$ strongly converging to $u_{\circ}$ in $L^2(\O)$
and let $\{u^{m}\}$ be a sequence of solutions of system
(\ref{PFl}) corresponding to the initial data $\{u_\circ^m(x)\}$.
The existence and regularity of such solutions has been gained in
Proposition\,\ref{existencenonu}. In particular we have
\footnote{Note that we also have, as in the proof of Proposition\,\ref{existencenonu}, $$ \|u^m_t\|_{L^{\infty}
(0,T;L^2(\O))}+\|u^m _t\|_{L^{2}(0,T;W^{1,p}(\O))}+\|u^m_t\|_{L^{p'}(0,T;V')}
 \leq M_1(\mu,T,\O, \|u_\circ^m\|_{2,2})\,,
$$ but this
estimate is not uniform in $m$ and $\mu$. }
\be\label{fea}
\|u^m \|_{L^{\infty}(0,T;L^2(\O))}+\|u^m \|_{L^{p}(0,T;V)}
\leq M(T,\O, \|u_\circ\|_{2})\,,\ \forall m\in \N.\ee Moreover, from
 Proposition\,\ref{existpeso1}, \be\ba{ll}\label{feb}\vspace{1ex} \dy \|t^\frac 1p \nabla
u^m\|_{L^{\infty}(0,T;L^p(\O))}+\|t\,u^m_t\|_{L^{\infty}(0,T;L^2(\O))}+\|t^\frac{p+2}{2p}\,\nabla
u^m_t\|_{L^{2}(0,T; L^p(\O))} \\ \dy\hskip7.5cm \leq M_2(T,\O, \|u_\circ\|_{2})\,,\ \forall m\in \N,\ea\ee and,
for $p>p_\circ$,

\be\label{fec1}\ba{ll}\vs1\dy t^\alpha \|\nabla u^m(t)\|_2^2\leq C\,
\mu^\frac{n(2-p)^2}{2(2p-n(2-p))}\,B^\frac 2p(\mu, u_\circ)\,
t^{\frac{n(2-p)^2}{p(2p-n(2-p))}}\\ \hskip2cm \dy+C\,(B(\mu,
u_\circ))^\frac{4-n(2-p)}{2p-n(2-p)} +c B(\mu, u_\circ)
t^\alpha(1+t)\,,\ea\ee
\be\label{fec2}\ba{ll}\dy\int_0^t\tau^{\alpha+\frac{2-p}{p}} \|D^2
u^m(\tau)\|_p^2\, d\tau
 \leq
C(B(\mu, u_\circ),T).\ea\ee

As the estimates in the norms \eqref{fea}--\eqref{fec2} are uniform in $m$, we
can extract a subsequence, still denoted
 by $\{u^m\}$, weakly or weakly-* converging in the same norms. Further,
 $u^m\in L^p(0,T; V)$, implies $a(\mu, u^m)\,\nabla u^m \in
L^{p'}(0,T; L^{p'}(\O))$, hence $A_\mu(u^m) \in L^{p'}(0,T; V')$
and  $u^m_t\in L^{p'}(0,T; V')$ and, as $m\to \infty$,
$$A_\mu(u^m)\rightharpoonup \overline\chi\ \textrm{ in } L^{p'}(0,T;
V')\ \textrm{weakly}.$$ Moreover,
we have
\be\label{strong1}\ba{ll}\dy
\|u^m(t)-u^{k}(t)\|_2^2+
\int_0^t\!(A_\mu(u^m)-A_\mu(u^{k}),
u^m-u^{k})d\tau \\ \dy \hskip6cm
\leq
\|u^{m}_\circ-u^{k}_\circ\|_2^2,
\forall t\geq 0, \forall m,k\in
\N.\ea\ee
 From this inequality and taking into account the $L^2$-strong
convergence of the sequence $\{u_\circ^m\}$ to $u_\circ$ and the monotonicity of the
operator, it follows the strong convergence of the sequence $\{u^m\}$ to
$u$ in $L^2(\O)$, uniformly in $t\geq 0$. Now, exactly as in
\cite{lions} one proves that the limit $u$ is the unique solution of
\eqref{PFl}, corresponding to the initial data $u_\circ\in L^2(\O)$.
\par If $u_\circ\in W_0^{1,2}(\O)$, we can reason as before, choosing a sequence $\{u_\circ^m\}$ strongly converging to $u_\circ$ in $W_0^{1,2}(\O)$, and replacing estimates \eqref{fec1}, \eqref{fec2} by
the following one, which, from Proposition\,\ref{existpeso1}, holds for any $p>\frac 32$
\be\label{ago22}\|\, \nabla u^m\|_{C(0,T;L^{2}(\O))}+
\|\, u^m\|_{L^{2}(0,T;W^{2,p}(\O))}\leq M_3(T,\O, \|u_\circ\|_{1,2})\,.\ee
We omit further details.
\vskip0.3cm
{\bf Proof of Theorem\,\ref{existencenomu}} -  From
Proposition\,\ref{DDAAGG}, for any fixed
$\mu>0$, there exists a unique solution of \eqref{PFl}. Let us denote by
$\{u^\mu\}$ the sequence of solutions of \eqref{PFl} for the different
values of $\mu>0$.    This sequence satisfies the bounds in \eqref{fea}--\eqref{fec2},
 uniformly in $\mu$. Hence, we
can extract a subsequence, still denoted
 by $\{u^\mu\}$, weakly or weakly-$*$ converging in the same norms,
 as $\mu$ goes to zero. In particular, in the limit as $\mu$ goes to zero,  estimate \eqref{fec1}
 gives
 $$ \ba{ll}\vs1\dy t^\alpha \|\nabla
u(t)\|_2^2\leq C\,\|u_\circ\|_2^\frac{2(4-n(2-p))}{2p-n(2-p)} +c
\|u_\circ\|_2^2 t^\alpha(1+t)\,,\ea$$ hence $t^\frac{\alpha}{2}
\nabla u\in C(0,T; L^2(\O))$.
 Let us show that the limit $u$ is the unique solution of system \eqref{PF}.
Recall that, from \eqref{oper1}, $$A(\psi):=\,-\nabla
\cdot \left( |\nabla \psi|^\frac{p-2}{2}\,\nabla \psi\right)\,. $$
Since
$u^\mu$ belongs to $L^p(0,T; V)$ and it is $\mu$-uniformly bounded, then $a(\mu, u^\mu)\,\nabla u^\mu \in
L^{p'}(0,T; L^{p'}(\O))$, and therefore $A_\mu(u^\mu) \in L^{p'}(0,T; V')$ and, along a subsequence,

$$A_\mu(u^\mu)\rightharpoonup \chi\ \textrm{in } L^{p'}(0,T; V')\
\textrm{weakly}.$$ Let us show that $\chi=A(u)$. Firstly we observe
that $$ |\,(\mu+|\nabla \psi|^2)^\frac{p-2}{2}\,\nabla \psi-|\nabla
\psi|^{p-2}\nabla\psi\,|^{p'} \to 0, \ \mbox{\ a.e. in }\
(0,T)\times \O,$$ and
$$|\,(\mu+|\nabla \psi|^2)^\frac{p-2}{2}\,\nabla \psi-|\nabla
 \psi|^{p-2}\nabla\psi\,|^{p'}\leq 2^p|\nabla \psi|^p. $$
The Lebesgue dominated convergence theorem ensures that
$$
(\mu+|\nabla \psi|^2)^\frac{p-2}{2}\,\nabla \psi \to
|\nabla \psi|^{p-2}\nabla \psi\ \mbox{ in } L^{p'}((0,T)\times \O) \ \mbox{ strongly} \,,$$
hence  \be\label{strongcA}A_\mu(\psi)\to A(\psi) \mbox{ in } L^{p'}(0,T; V') \ \mbox{ strongly} \,.\ee
  Set
$$X_{\mu}:=\int_0^t\, (A_\mu(u^\mu)-A_\mu(\psi), u^\mu-\psi)\,d\tau\,\geq 0,\ \forall
\psi\,\in L^p(0,T,V)\,.$$
By using that $u^\mu$ is a solution of
\eqref{PFl},  we write $X_\mu$ as follows
$$\ba{ll}\vs1\dy X_{\mu}\!\!\!\!&= \dy \int_0^t\ (A_\mu(u^\mu),u^\mu)-(A_\mu(u^\mu),\psi)-(A_\mu(\psi),
u^\mu-\psi)\, d\tau \\\vs1\dy &=\dy \frac 12 \|u_\circ\|_2^2-\frac 12 \|u^\mu(t)\|_2^2-\int_0^t\ [ (A_\mu(u^\mu),\psi)+(A_\mu(\psi),
u^\mu)-(A_\mu(\psi),\psi)]\, d\tau\,.\dy \ea
$$
Let us pass to the lim inf. Observe that, since
$$(A_\mu(\psi),u^\mu)- (A(\psi),u) = (A_\mu(\psi)-A(\psi),
u^\mu)+(A(\psi),u^\mu-u),$$
from the strong convergence \eqref{strongcA}, the uniform bound
of $u^\mu$ in $L^p(0,T;V)$ and the weak convergence of $ u^\mu$
to $u$ in $L^p(0,T;V)$, we get
$$\lim_{\mu\to 0} \int_0^t  (A_\mu(\psi),u^\mu) \, d\tau=
 \int_0^t  (A(\psi),u)\,d\tau \,.$$ Hence
$$\ba{ll}\dy 0\leq \liminf X_{\mu}\leq \dy\ \frac 12 \|u_\circ\|_2^2
-\frac 12 \|u(t)\|_2^2\dy - \int_0^t (\chi,\psi) +(A(\psi),
u-\psi)d\tau .
\ea
$$
On the other hand, it is easy to see that
the limit $u$ satisfies \be\label{appp2mu}\int_0^t(\chi, u)\, d\tau
=\,\frac 12 \|u_\circ\|_2^2-\frac 12 \|u(t)\|_2^2\,.\ee
Therefore we obtain
$$\int_0^t (\chi-A(\psi), u-\psi)\,d\tau \geq 0, \  \ \forall\psi\in L^p(0,T; V)\,.$$
Taking $\psi=u+\lambda w$, for $\lambda>0$ and for some $w\in
L^{p}(0,T, V)$, and then letting $\lambda$ tend to zero the thesis
follows.\chiu
\vskip0.3cm
{\bf Proof of Corollary\,\ref{corexistencenomu}} - From Proposition\,\ref{DDAAGG}, for any fixed
$\mu>0$, there exists a solution of \eqref{PFl}.  Let us denote by
$\{u^\mu\}$ the sequence of solutions of \eqref{PFl} for the different
values of $\mu>0$.    This sequence satisfies the bounds in \eqref{fea}, \eqref{feb} and \eqref{ago22},
 uniformly in $\mu$. Hence, we
can extract a subsequence, still denoted
 by $\{u^\mu\}$, weakly or weakly-* converging in the same norms,
 as $\mu$ goes to zero.
  That the limit $u$ is the unique solution of system \eqref{PF} can be proved as in the proof of Theorem\,\ref{existencenomu}.

\section*{\normalsize 5. A crucial estimate: $u_t\in L^{\infty}(0,T; L^q(\O))$}
\renewcommand{\theequation}{5.\arabic{equation}}
\renewcommand{\thetho}{5.\arabic{tho}}
\renewcommand{\thedefi}{5.\arabic{defi}}
\renewcommand{\therem}{5.\arabic{rem}}
\renewcommand{\theprop}{5.\arabic{prop}}
\renewcommand{\thelemma}{5.\arabic{lemma}}
\renewcommand{\thecoro}{5.\arabic{coro}}
\setcounter{equation}{0} \setcounter{coro}{0} \setcounter{lemma}{0}
\setcounter{defi}{0} \setcounter{prop}{0} \setcounter{rem}{0}
\setcounter{tho}{0}
\begin{prop}\label{ut1} {\sl Let $p>p_\circ$, with $p_\circ$ given in \eqref{restrp},
 and $q\in  [2,\frac{7-3p}{4-2p})$.
 Let $u$ be the unique solution of \eqref{PF} corresponding to
$u_\circ\in L^2(\O)$. Then $t^{1+\gamma}\,u_t\in
L^\infty(0,T;L^{q}(\O))$,
 with  $\gamma=\gamma(q')$ given by \eqref{Maxpostraa}. Moreover the following estimate holds
 \be\label{estvt}\|u_t(t)\|_{q}\leq \frac{c}{t}{\null_{1+\gamma}}\,
\|u_\circ\|_2^{(2-p)\gamma+1},\, \forall t\in (0,T).\ee}
\end{prop}
\Pr
Let $\dy\{ u_\circ^m(x)\}$ be a sequence
in $C_0^\infty(\O)$ strongly converging to $u_{\circ}$ in $L^2(\O)$
and let $\{v^{m}\}$ be the sequence of solutions of system
(\ref{PFepv}) corresponding to the initial data $\{u_\circ^m(x)\}$.
Under our assumptions on $p$, from Proposition\,\ref{existence}
 and Corollary\,\ref{pregcor} we know that,
 for any data $u_\circ^m$, there exists a unique solution $v^m$
 of \eqref{PFepv} such that $t^{\beta_2(p)}\,v^n
 \in L^2(0,T; W^{2,p}(\O))$ and $t\,v^m_{t} \in L^{\infty}(0,T; L^{2}(\O))$.
 Therefore $v^m$ satisfies system \eqref{PFepv} a.e. in $(0,T)\times
 \O$. In the sequel, for simplicity we suppress the superscript $m$.\par
Let us mollify equation \eqref{PFepv} with respect to $t$, and
denote the mollifier by $J_\delta$. We can derive the regularized
system \eqref{PFepv} with respect to $t$ and obtain, a.e.
 in $\O\times (0,T)$,
\be\label{vd} \pa_t\, J_\delta (v_{t})-\nu \pa_t \Delta (J_\delta\, v)-\nabla\cdot
\pa_t\,J_\delta\left((\mu+|\nabla
v|^2)^\frac{(p-2)}{2}\nabla v\right)= 0. \ee For any
fixed $\eta$, let  $\vp^\eta(s,x)$ be the unique solution of system
\eqref{AD1} with $\nu>0$, corresponding to an initial data
$\vp_\circ \in C_0^\infty(\O)$, and $B_\eta$, given by
\eqref{beeta}, with $v$ as above. For simplicity we also drop the
superscript $\eta$.
  From Lemma\,\ref{newag1}, $\vp\in L^{2}(0,t; W_0^{1,2}(\O))$ and $\vp_s\in
L^{2}((0,t)\times \O)$. Therefore, setting $\tau=t-s$, $s\in
(0,\frac t2)$, we can multiply equation \eqref{vd} by $\vp(t-\tau)$,
for $\tau\in (\frac t2,t)$, and integrate in $(\frac t2, t)\times
\O$. Integrations by parts give the following identity
\be\label{deltato}\ba{ll}\vs1\dy (J_\delta
(v_t)(t),\vp(0))=(J_\delta (v_t)(\frac t2),\vp(\frac t2))+
\int_{\frac t2}^t(J_\delta (v_{\tau})(\tau),\vp_\tau(t-\tau))\\
\vs1\dy +\nu\!\int_{\frac t2}^t\!(J_\delta(\nabla
v_{\tau})(\tau),\nabla \vp(t-\tau))d\tau \!+\! \int_{\frac t2}^t\!\!
\left(\!J_\delta\!\left(\frac{\nabla v_\tau}{(\mu+|\nabla
v|^2)}{\null_\frac{(2-p)}{2}}\!\!\right)\!(\tau), \nabla \vp(t-\tau)\!\!\right)\!d\tau\\
\dy+ (p-2)\int_{\frac t2}^t \left(J_\delta \left(\frac{(\nabla
v\otimes\nabla v)\cdot\nabla v_\tau}{(\mu+|\nabla
v|^2)^\frac{(4-p)}{2}}\right)(\tau),\nabla
\vp(t-\tau)\right)d\tau.\ea\ee Let us pass to the limit as $\delta$
goes to zero. Note that from Corollary\,\ref{existencecor1} one gets
$t\,v_t\in C([0,T);L^2(\O))$. Hence, since, from Lemma\,\ref{adgt},
we also have $\vp\in C([0,t); L^2(\O))$, we get
$$\lim_{\delta\to 0} (J_\delta (v_t)(\frac t2),\vp(\frac t2))= (v_t(\frac t2),\vp(\frac t2)),$$
as well as
$$\lim_{\delta\to 0} (J_\delta (v_t)(t),\vp_\circ)= (v_t(t),\vp_\circ).$$
Further, from  Corollary\,\ref{existencecor1}, one has $\tau^\frac 12
\,v_\tau\in L^2((0,T)\times \O)$, hence $J_\delta (t^\frac 12\,v_t)$
strongly converges to $t^\frac 12 \,v_t$ in $L^2((0,T)\times \O)$, as $\delta\to 0$.
Since,
 from Lemma\,\ref{newag1}, one has $\tau^{-\frac 12}\,\vp_\tau(t-\tau)\in L^2((\frac t2,t)\times \O)$, we get
 $$\lim_{\delta\to 0}
 \int_{\frac t2}^t(J_\delta
(v_{\tau})(\tau),\vp_\tau(t-\tau))d\tau =
 \int_{\frac t2}^t(v_{\tau}(\tau),\vp_\tau(t-\tau))d\tau\,.$$ Recall that, from Corollary\,\ref{existencecor1},
$$\|J_\delta (\nabla
v_\tau)\|_{L^2(\frac t2,T;L^2(\O))}\leq \|\nabla v_\tau\|_{L^2(\frac
t2,T;L^2(\O))}\leq \frac{c}{t}B^\frac 12(\mu, u_\circ^m).$$Therefore
 $$\|J_\delta \frac{\nabla
v_\tau}{(\mu+|\nabla v|^2)}{\null_\frac{(2-p)}{2}}\|_{L^2((\frac
t2,T)\times\O)}\leq \|\frac{\nabla v_\tau}{(\mu+|\nabla
v|^2)}{\null_\frac{(2-p)}{2}} \|_{L^2((\frac t2,T)\times\O)}\leq
\frac{c}{t}\mu^{p-2}\,B^\frac 12(\mu, u_\circ^m),$$ and, similarly,
 $$\|J_\delta \frac{(\nabla v\otimes\nabla
v)\cdot (\,\nabla v_\tau)}{(\mu+|\nabla
v|^2)^\frac{(4-p)}{2}}\|_{L^2((\frac t2,T)\times\O)}\!\leq
\|\frac{\nabla v_\tau}{(\mu+|\nabla
v|^2)}{\null_\frac{(2-p)}{2}}\|_{L^2((\frac t2,T)\times\O)}\leq
\frac {c}{t}\mu^{p-2}B^\frac 12(\mu, u_\circ^m).$$ Moreover,
observing   that $\tau^{-1}\,\nabla \vp_\tau\in L^2(\frac t2,t;
L^2(\O))$, we can pass to the limit as $\delta$ tends to zero in the
last two integrals on the righ-hand side of \eqref{deltato}.
Summarizing, in the limit  as $\delta$ goes to zero, \eqref{deltato}
gives \be\label{deltato1}\ba{ll}\dy (v_t(t),\vp(0))=(v_t(\frac
t2),\vp(\frac t2))+\!\int_{\frac
t2}^t\!(v_{\tau}(\tau),\vp_\tau(t-\tau))d\tau\\
\vs1\dy\hskip0.5cm +\nu\!\int_{\frac t2}^t\!(\nabla
v_{\tau}(\tau),\nabla \vp(t-\tau))d\tau+ \int_{\frac t2}^t
\!\left(\frac{\nabla v_\tau(\tau)}{(\mu+|\nabla
v(\tau)|^2)}{\null_\frac{(2-p)}{2}}, \nabla
\vp(t-\tau)\right)d\tau\\\hskip3cm \dy \ + (p-2)\int_{\frac t2}^t
\!\!\left(\frac{(\nabla v(\tau)\otimes\nabla v(\tau))\cdot\nabla
v_\tau(\tau)}{(\mu+|\nabla v(\tau)|^2)^\frac{(4-p)}{2}},\nabla
\vp(t-\tau)\right)d\tau.\ea\ee Let us  write the last two integrals
on the right-hand side of \eqref{deltato1} as follows
$$\ba{ll}\dy \vs1 \int_{\frac t2}^t \!\!\left(\frac{\nabla
v_\tau(\tau)}{(\mu+|\nabla v|^2)^\frac{(2-p)}{2}}, \nabla
\vp(t-\tau)\right)d\tau+ (p-2)\!\int_{\frac t2}^t
\!\left(\frac{(\nabla v\otimes\nabla v)\cdot\nabla
v_\tau}{(\mu+|\nabla v|^2)^\frac{(4-p)}{2}},\nabla
\vp(t-\tau)\right)d\tau \\\vs1 \dy = \int_{\frac t2}^t \left( \nabla
v_\tau(\tau), B_\eta((\tau,x)\vp(t-\tau)\right)d\tau \\\quad\vs1 \dy
+
 \int_{\frac t2}^t \left( \nabla v_\tau (\tau),\left[\frac{1}{(\mu+|\nabla v|^2)^\frac{2-p}{2}}-\frac{1}{
(\mu+|J_\eta(\nabla v)|^2)^\frac{2-p}{2}}\right]
\nabla\vp(t-\tau)\right) d\tau\\\dy\quad +(p-2)\int_{\frac t2}^t
\left(\nabla v_\tau,  \left[\frac{(\nabla v)\otimes (\nabla
v)}{(\mu+|\nabla v|^2)^\frac{4-p}{2}}- \frac{J_\eta(\nabla v)\otimes
J_\eta(\nabla v)}{(\mu+|J_\eta(\nabla
v)|^2)^\frac{(4-p)}{2}}\right]\cdot\nabla \vp
(t-\tau)\right)d\tau\,. \ea
$$ Then, since $\vp$ is solution of \eqref{AD1} with $B_\eta$ as in \eqref{beeta},  identity
\eqref{deltato1} becomes \be\label{deltato2}\ba{ll}\dy\vs1
(v_t(t),\vp_\circ)=(v_t(\frac t2),\vp(\frac t2))\\\vs1
\dy \quad +
 \int_{\frac t2}^t \left(\nabla v_\tau (\tau),\left[\frac{1}{(\mu+|\nabla v|^2)}{\null_{\frac{2-p}{2}}}-\frac{1}{
(\mu+|J_\eta(\nabla v)|^2)}{\null_{\frac{2-p}{2}}}\right]
\nabla\vp(t-\tau)\right) d\tau\\\dy\quad\vs1 +(p-2)\!\!\!\int_{\frac
t2}^t\!\! \left(\!\nabla v_\tau,  \left[\frac{(\nabla v)\otimes
(\nabla v)}{(\mu+|\nabla v|^2)} {\null_{\frac{4-p}{2}}}-
\frac{J_\eta(\nabla v)\otimes J_\eta(\nabla
v)}{(\mu+|J_\eta(\nabla v)|^2)}{\null_{\frac{4-p}{2}}}\right]\!\cdot\nabla \vp (t-\tau)\!\right)\!d\tau\\
\dy \ = (v_t(\frac t2),\vp(\frac t2))+ I^1_\eta+ I^2_\eta .\ea\ee We
claim that the integrals $I^1_\eta$ and $I^2_\eta$ go to zero along
a subsequence. This follows from Lemma\,\ref{compattificazione}, with
$h^\eta=\tau^{-1}\,\nabla \vp^\eta $, which is bounded in
$L^2((\frac t2,t)\times \O)$, uniformly in $\eta>0$, and
$\psi=\tau\, \nabla v_\tau$ which is also in $L^2((\frac t2,T)\times
\O)$, due to Corollary\,\ref{existencecor1}.
\par
Therefore,
 by using estimate \eqref{Maxpostaa} in  Lemma\,\ref{newag2}
  with $r=q'$, $s=\frac t2$ and $M=B$, and using that, from {\it i)}
 in Corollary\,\ref{existencecor1},
  $tv_t\in L^\infty(0,T; L^2(\O))$(see estimate
 \eqref{a16}), then passing to the limit as $\eta$ goes to zero in \eqref{deltato2}
  we get $$\ba{ll}\dy (v_t(t),\vp_\circ)=(v_t(\frac
t2),\vp(\frac t2))\\ \vs1\dy\hskip3cm\leq \|v_t(\frac t2)\|\,
\|\vp(\frac t2)\|\leq \frac{c}{t^{1+\gamma(q')}}\|\vp_\circ\|_{q'}\,
\left( \|u_\circ^m\|_2^2+c(\O,T)\mu^\frac
p2\right)^{\!\frac{(2-p)\gamma+1}{2}},\ea$$ with $\gamma(q')$ given in
\eqref{Maxpostraa}. Let us give back the superscript $m$ to the
sequence $\{v^m\}$. The last estimate ensures that $t^{1+\gamma}\,
v^m_t\in L^\infty(0,T;L^{q}(\O))$ and that the following estimate
holds
 \be\label{estvtv}\|v^m_t(t)\|_{q}\!\leq \frac{c}{t}{\null_{1+\gamma(q')}}\!
\left(\|u_\circ\|_2^2+c(\O,T)\mu^\frac
p2\!\right)^{\!\!\frac{(2-p)\gamma(q')+1}{2}}\!\!\!\!,\mbox{ uniformly in
}\nu>0, \mu>0, m\in \N.\ee
 Now we pass to the limit, firstly as $\nu$ goes to zero, thus
 obtaining the  sequence of solutions $\{u^m\}$ of Proposition\,\ref{existencenonu},
 whose time derivative also satisfies
  estimate \eqref{estvtv}. Then we pass to the limit as $m$
  tends to infinity, finding the solution $u=u^\mu$ of system
   \eqref{PFl} and satisfying \eqref{estvtv}. Finally, we pass
   to the limit as $\mu$ goes to zero, thus obtaining that the
   solution $u$, found in Theorem\,\ref{existencenomu},  satisfies
   \eqref{estvt}.\chiu
\begin{rem} {\rm Note that the main reason that leads us to resort
 to system \eqref{PFepv} as approximation of system \eqref{PFl} is
 due to the possibility of  handling  the integrals $I^1_\eta$ and
 $I^2_\eta$.
}\end{rem}

\section*{\normalsize 6. Proof of Theorem\,\ref{regularitynomu},
 Theorem\,\ref{regularitynomun2}
  and their corollaries}
\renewcommand{\theequation}{6.\arabic{equation}}
\renewcommand{\thetho}{6.\arabic{tho}}
\renewcommand{\thedefi}{6.\arabic{defi}}
\renewcommand{\therem}{6.\arabic{rem}}
\renewcommand{\theprop}{6.\arabic{prop}}
\renewcommand{\thelemma}{6.\arabic{lemma}}
\renewcommand{\thecoro}{6.\arabic{coro}}
\setcounter{equation}{0} \setcounter{coro}{0} \setcounter{lemma}{0}
\setcounter{defi}{0} \setcounter{prop}{0} \setcounter{rem}{0}
\setcounter{tho}{0} In order to prove Theorem\,\ref{regularitynomu} and its Corollary\,\ref{corregularitynomu},
firstly we recall a regularity result obtained in
\cite{CMellittico}, related to the singular elliptic
$p$-Laplacean system \be\label{real}\ba{ll}\vs1\dy-\nabla \cdot
\left(\,|\nabla u| ^{p-2}\,\nabla u\,\right)= f\,,&\dy \mbox{ in }
\Omega\,,\\\dy\hskip2.8cm u=0 \,, &\dy \mbox{ on } \partial \O\,,\ea \ee with
$p\in (1,2)$. We recall that, if $v \in W^{2,2}(\O)\cap
W_0^{1,2}(\O)$, then the following estimate holds
 \be\label{C3rich}\|D^2 v\|_2\, \leq
 H \|\Delta
v\|_2,\ee where the constant $H$ depends on the size of $\O$. If
$\O$ is a convex domain the inequality holds with $H=1$. For
details we refer to to \cite{L} or \cite{LU}. We recall that
\be\label{ovp}\overline p=2-\mbox{\large$\frac1H\,,$}\ee where $H$ is the above
constant.
The results in \cite{CMellittico} can be stated as follows.

\begin{tho}\label{threal}
{\sl Let $\O\subset \R^n$, $n\geq 3$, be convex.
Assume that $f\in L^q(\O)$, with $q\in [ \frac{2n}{n(p-1)+2(2-p)}, \infty)$. There exists a
unique solution $u\in
W_0^{1,2}(\O)\cap
 W^{2,\widehat q}(\O)$ of system \eqref{real}, with
\begin{equation}\label{mainestu}\|\, u \,\|_{2,\widehat q}
\leq C \|f\|_{q}^\frac{1}{p-1}\,,
\end{equation} where $C$ is a
constant independent of $u$ and
$\widehat q=\widehat q(q)= \frac{nq(p-1)}{n-q(2-p)}$ if  $q<\,n$,
$\widehat q <n$ if $q =\, n$, and $\widehat q=q$ if $q>\, n$.
 \par The same result holds for
non-convex domains $\O$, provided
that $p>\overline p$}.
\end{tho}

\begin{tho}\label{threaln2}
{\sl Let  $\O\subset \R^2$, be convex.
Assume that $f\in L^q(\O)$, with $q>2$.
There exists a
unique solution $u\in
W_0^{1,2}(\O)\cap
 W^{2,q}(\O)$ of system \eqref{real}, with
\begin{equation}\label{mainestun2}\|\, u \,\|_{2,q}
\leq C \|f\|_{q}^\frac{1}{p-1}\,,
\end{equation} where $C$ is a
constant independent of $u$.
\par The same result holds for
non-convex domains $\O$, provided
that $p>\overline p$}.
\end{tho}

{\bf Proof of Theorem\,\ref{regularitynomu}} -
Since $p>p_\circ$, we can apply Theorem\,\ref{existencenomu}, and find that the unique solution
 of system \eqref{PF} satisfies the following system, a.e. in
$t\in (0,T)$,
$$\ba{ll}\vs1\dy-\nabla \cdot \left(\,|\nabla u| ^{p-2}\,\nabla
u\,\right)= u_t\,,&\dy \mbox{ in } \Omega\,,\\\dy\hskip2.8cm u=0 \,,
&\dy \mbox{ on } \partial \O\,.\ea $$ We set \be\label{ovq} \ov q:=\mbox{\large$
\frac{2n}{n(p-1)+2(2-p)}\,,$}\ee and observe that $\ov q\geq 2$ and,
since $p>p_1$,
then $\ov q\leq \frac{7-3p}{4-2p}$. From Proposition\,\ref{ut1}, for
$p>p_\circ$ and $q\in [\,\ov q,\frac{7-3p}{4-2p}\,]$, we have
$t^{1+\gamma} u_t\in L^\infty(0,T; L^q(\O))$, where
$\gamma=\frac{n(q-2)}{q(2p-2n+np)}$,  and $u_t$ satisfies estimate
\eqref{estvt}. By applying the above Theorem\,\ref{threal} we find that: if $\O$ is convex, then $u\in
W^{2,\widehat q}(\O)$, a.e. in $t\in (0,T)$; if $\O$ is not convex,
the same result holds if $p>\ov p$ too. In both cases,
using estimate \eqref{estvt} in \eqref{mainestu} we find
\be\label{percor1}\|u \|_{2,\widehat q} \leq C \,\|u_t\|_{q}^\frac{1}{p-1}\leq
\frac{c}{t}{\null_{1+\gamma}}\,
\|u_\circ\|_2^{\frac{(2-p)\gamma+1}{p-1}},\, \mbox{ a.e. in } t\in
(0,T)\,.
\ee
\chiu
\vskip0.1cm
For the proof of Theorem\,\ref{regularitynomun2} we argue exactly in the same way,
employing Theorem\,\ref{threaln2}
 in place of Theorem\,\ref{threal}.
For the sake of completeness we give the proof.
\vskip0.1cm{\bf Proof of Theorem\,\ref{regularitynomun2}} -
From Theorem\,\ref{existencenomu}, we can write the following system, a.e. in
$t\in (0,T)$,
$$\ba{ll}\vs1\dy-\nabla \cdot \left(\,|\nabla u| ^{p-2}\,\nabla
u\,\right)= u_t\,,&\dy \mbox{ in } \Omega\,,\\\dy\hskip2.8cm u=0 \,,
&\dy \mbox{ on } \partial \O\,.\ea $$ From Proposition\,\ref{ut1}, for
$p>p_\circ$ and $q\in(2,\frac{7-3p}{4-2p}]$, we have
$t^{1+\gamma} u_t\in L^\infty(0,T; L^q(\O))$, where
$\gamma=\frac{(q-2)}{2q(p-1)}$,  and $u_t$ satisfies estimate
\eqref{estvt}. We apply Theorem \ref{threaln2} and we find that: if $\O$
is convex, then $u\in W^{2,\widehat q}(\O)$, with
$\widehat q=q$, a.e. in $t\in (0,T)$, ; if $\O$ is not convex, the
same result holds for $p>\ov p$ too. In both cases,
using estimate \eqref{estvt} in
\eqref{mainestun2}, we find
\be\label{percor2}\|u \|_{2,\widehat q} \leq C \,\|u_t\|_{q}^\frac{1}{p-1}\leq
\frac{c}{t}{\null_{1+\gamma}}\,
\|u_\circ\|_2^{\frac{(2-p)\gamma+1}{p-1}},\, \mbox{ a.e. in } t\in
(0,T)\,.
\ee\chiu\vskip0.2cm

{\bf Proof of Corollary\,\ref{corregularitynomu}} -
Firstly we prove estimate \eqref{holdertxa}. Since $q_\circ\in(\frac
np,n]$, then Theorem\,\ref{regularitynomu} ensures that
$u\in
L^\infty(\ve,T;W^{2,\widehat
q_\circ}(\O))$ and $u_t\in
L^\infty(\ve,T;L^{q_\circ}(\O)),$ with
$\widehat
q_\circ=\frac{nq_\circ(p-1)}{n-q_\circ(2-p)}\in (\frac n2,n)$ if
$q_\circ\in (\frac np,n)$, $\widehat q_\circ\in (\frac n2,n)$ if $q_\circ=n$.  From Lemma\,\ref{solonnl} with $m=0$ and $\ov q=\widehat q_\circ$
we obtain
$$\big[\,u\,\big]_{\lambda_\circ,t,x}\leq C \left[\sup_{(\ve,T)} (\|u_t(t)\|_{\widehat q_\circ}+ \|D^2 u(t)\|_{\widehat q_\circ})
+\sup_{(\ve,T)}\|u(t)\|_{\widehat q_\circ}\right], $$
$\lambda_\circ=2-\frac
{n}{\widehat q_\circ}$,
from which estimate \eqref{holdertxa} easily follows by increasing the right-hand side via estimate \eqref{percor1}.
Assume now that $p>\max\{p_\circ,
 \frac{4n-7}{2n-3}\}$ and $q_1\in (n,
 \frac{7-3p}{4-2p}]$.
 Under the assumption $p>\frac{4n-7}{2n-3}$,
$\frac{7-3p}{4-2p}>n$ is ensured.
 Therefore we can apply
Theorem\,\ref{regularitynomu} with $q_1=q>n$, and obtain
$u_t\in L^\infty(0,T;L^{q_1}(\O))$, $u\in L^\infty(0,T;W^{2,q_1}(\O))$. By applying Lemma\,\ref{solonnl} with $m=1$ and $\ov q=q_1$  we get estimate \eqref{holdertx}. \chiu\vskip0.2cm
{\bf Proof of Corollary\,\ref{corregularitynomu2}}  -
Theorem\,\ref{regularitynomun2}
gives, for $q>2$,   $u\in L^\infty(0,T;W^{2,q}(\O))$ and $u_t\in L^\infty(0,T;L^{q}(\O))$. Hence, by applying Lemma\,\ref{solonnl} with $m=1$ and $\ov q=q$  we get estimate \eqref{holdertxb}.

 \chiu\vskip0.2cm
 \section*{\normalsize 7. The Maximum Modulus Theorem: Proof of Theorem\,\ref{mmtp}}
\renewcommand{\theequation}{7.\arabic{equation}}
\renewcommand{\thetho}{7.\arabic{tho}}
\renewcommand{\thedefi}{7.\arabic{defi}}
\renewcommand{\therem}{7.\arabic{rem}}
\renewcommand{\theprop}{7.\arabic{prop}}
\renewcommand{\thelemma}{7.\arabic{lemma}}
\renewcommand{\thecoro}{7.\arabic{coro}}
\setcounter{equation}{0} \setcounter{coro}{0} \setcounter{lemma}{0}
\setcounter{defi}{0} \setcounter{prop}{0} \setcounter{rem}{0}
\setcounter{tho}{0}

\begin{prop}\label{ulinfty} {\sl Let  $\mu> 0$.
Let $u$ be the solution of \eqref{PFl} corresponding to
$u_\circ\in L^\infty(\O)$. Then $$\|u(t)\|_\infty\leq
\|u_\circ\|_\infty.$$}
\end{prop}
\Pr We argue as in Proposition\,\ref{ut1}. Hence we consider
a sequence  $\dy\{ u_\circ^m(x)\}\in C_0^\infty(\O)$ strongly
convergent to $u_{\circ}$ in $L^2(\O)$, and such that
$\|u_\circ^m\|_\infty\leq \|u_\circ\|_\infty$, the
 sequence  $\{v^{m,\nu}\}$  of solutions of system
(\ref{PFepv}) corresponding to the initial data $\{u_\circ^m\}$, for
which Proposition\,\ref{existence} and Corollary\,\ref{existencecor1}
hold.
 In the sequel, for simplicity we suppress the superscripts $m$
 and $\nu$ for the sequence $\{v^{m,\nu}\}$.
Let us consider, for a fixed $\eta>0$,  the solution $\vp^\eta(s,x)$
of system \eqref{AD1} corresponding to a data $\vp_\circ\in
C_0^\infty(\O)$, where $B_\eta$ is given by \eqref{beeta} with $b=0$,
 and set
$\hat\vp^\eta(\tau)=\vp^\eta(t-\tau)$. Note that, from estimate \eqref{os4},
 $v$ satisfies \eqref{beeta2}. Then, using $\hat\vp^\eta(\tau)$ as test
function in the weak formulation of \eqref{PFepv}, we
have
 \be\label{a8}\ba{ll}\dy\vs1 (v(t),\vp_\circ)-(v(0),\vp^\eta(t))
-\int_0^t (v(\tau),\hat\vp_\tau^\eta(\tau))d\tau
\\\dy\hskip2cm
=-\nu\int_0^t(\nabla v(\tau),\nabla
\hat\vp^\eta(\tau))d\tau-\int_0^t(a(\mu,v(\tau)) \nabla
v(\tau),\nabla \hat\vp^\eta(\tau))d\tau .\ea\ee Writing the second
term on the right-hand side of \eqref{a8} as
$$\ba{ll}\dy \vs1\int_0^t(a(\mu,v)
\nabla v(\tau),\nabla \hat\vp^\eta(\tau))d\tau=\!\!&\dy\! \int_0^t(
\nabla v(\tau),a_\eta(\mu,v(\tau))\nabla
\hat\vp^\eta(\tau))d\tau\\&\dy +\int_0^t( \nabla v(\tau),\nabla
\hat\vp^\eta(\tau))[a(\mu,v) -a_\eta(\mu, v)]\,d\tau\,,\ea
$$ and recalling that $\vp^\eta$ is solution of \eqref{AD1}, identity
\eqref{a8} becomes \be\label{a9}\ba{ll}\dy\vs1
(v(t),\vp(0))=(v(0),\vp^\eta(t))\\ \dy \hskip1cm -\int_0^t( \nabla
v(\tau),\nabla \hat\vp^\eta(\tau))[a(\mu, v(\tau))-a_\eta(\mu,
v(\tau))]\,d\tau=(v(0),\vp^\eta(t))+I_\eta .\ea\ee The integral
$I_\eta$ goes to zero, as $\eta$ goes to zero, along a subsequence,
thanks to Lemma\,\ref{compattificazione} with $h^\eta=\vp^\eta\in
L^2(0,t;L^2(\O))$, due to Lemma\,\ref{newag1} and with $\psi= \nabla v\in L^{2}(0,t; L^2(\O))$,
due to  Proposition\,\ref{existence}.
Finally, using \eqref{max11a} and then passing to the limit as $\eta$
tends to zero  in \eqref{a9}, along a subsequence, we get
$$|(v(t),\vp_\circ)| \leq \|u_\circ^m\|_{r'}\,\|\vp_\circ\|_{r}, \
\forall \vp_\circ\in C_0^\infty(\O).$$ Giving back the superscripts
to the sequence $\{v^{m,\nu}\}$,
this last estimate together with the bound $\|u^m_\circ\|\leq
\|u_\circ\|_\infty$ imply $$
\|v^{m,\nu}(t)\|_{r'}=\!\!\sup_{\vf_\circ\in C_0^\infty(\O) \atop
|\vf_\circ|_r=1}|(v^{m,\nu}(t),\vf_\circ)| \leq
\|u_\circ^m\|_{r'}\leq |\O|^\frac{1}{r'}\|u^m_\circ\|_\infty\leq
|\O|^\frac{1}{r'}\|u_\circ\|_\infty\,.$$ Since the the right-hand
side is uniform with respect to $r'$, letting $r'$ go to $\infty$,
we obtain $$ \|v^{m,\nu}(t)\|_{\infty}\leq \|u_\circ\|_{\infty}\,.$$
 Now we pass to the limit, firstly as $\nu$ goes to zero. We
 obtain the sequence of solutions $\{u^m\}$
  as in Proposition\,\ref{existencenonu}.
 From Proposition\,\ref{existence}, we have, uniformly in $\nu>0$,
 that $v^{m,\nu}\in L^\infty(0,T; V)$, $v^{m,\nu}_t\in L^2(0,T; L^p(\O))$.
 Since $V$ is compactly embedded in $L^s(\O)$, for any $p<s<\frac{np}{n-p}$,
 and $L^s(\O)$ is continuously embedded in $L^p(\O)$,  from Lemma\,\ref{aubin},
 the sequence $\{v^{m,\nu}\}$ converges to $u^m$, as $\nu$ goes to zero,
  strongly in $L^{s}((0,T)\times \O)$, hence almost everywhere in $t$, $x$,
   along to a subsequence. Therefore, along such a subsequence, we find
\be\label{abg0}\ba{ll}\dy \vs1 |u^m(t,x)|\leq |u^m(t,x)-v^{m,\nu}(t,x)|+
|v^{m,\nu}(t,x)|\\ \hskip5cm \leq |u^m(t,x)-v^{m,\nu}(t,x)|+\|v^{m,\nu}(t)\|_\infty,\ea \ee
 a.e. in $(t,x)\in (0,T)\times \O$. Passing to the limit on $\nu$, we easily get
 \be\label{abg1}
 \|u^m(t)\|_\infty\leq \|u_\circ\|_\infty\,.\ee
We now pass to the limit as $m$ tends to infinity. Then $\{u_\circ^m\}$
 strongly converges to $u_\circ$ and, from \eqref{strong1}, the corresponding
 sequence of solution converges strongly in $L^2(\O)$, uniformly in $t \geq 0$,
  hence, along a subsequence, a.e. in $(t,x)\in (0,T)\times \O$. Therefore,
  reasoning as in \eqref{abg0}, we find that the solution $u=u^\mu$ of system
   \eqref{PFl} satisfies \eqref{abg1}.
\chiu

\begin{prop}\label{ul2} {\sl Let $p>\frac{2n}{n+2}$ and $\mu> 0$.
 Let $u$ be the solution of \eqref{PFl} corresponding to $u_\circ\in L^q(\O)$, for some $q\in
[2,+\infty]$. Then
$$\|u(t)\|_\infty\leq c\,\left(\|u_\circ\|_2^2+ c(\O,T)\mu^\frac
p2\right)^{\!\frac{(2-p)\beta}{q}}
\|u_\circ\|_q\,t^{-\frac{2\beta}{q}},\quad \forall t\in (\ve,T),$$
with $\beta$ given in \eqref{Maxpostr}.}
\end{prop}

\Pr We consider a sequence  $\dy\{ u_\circ^m(x)\}\in C_0^\infty(\O)$
strongly converging to $u_{\circ}$ in $L^q(\O)$, the sequence
$\{v^{m,\nu}\}$  of solutions of system (\ref{PFepv}), corresponding
to the initial data $\{u_\circ^m\}$, satisfying Proposition\,\ref{existence} and Corollary\,\ref{existencecor1}.
Arguing as in
Proposition\,\ref{ulinfty}, we arrive at the following estimate
\be\label{aaa9}\dy\vs1
(v^{m,\nu}(t),\vp_\circ)=(v^{m,\nu}(0),\vp(t))+I_\eta .\ee From
Lemma\,\ref{newag2} with $r=q'$ and $M=B(\mu, u_\circ^m)$, and
introducing two arbitrary conjugate exponents $\ov s>2$ and $\ov
s'$,  we can estimate the first term on the right-hand side of
\eqref{aaa9} as
$$\ba{ll}\vs1 |(v^{m,\nu}(0),\vp(t))|\!\!\!&\dy\leq \|u_\circ^m\|_q\,
\|\vp(t)\|_{q'}\leq c\,(B(\mu, u_\circ^m))^\frac{(2-p)\beta}{q}
\!\|u_\circ^m\|_q\, \|\vp_\circ\|_1\, t^{-\frac{2\beta}{q}}\,,
\\& \dy \leq c\,|\O|^{\frac{1}{\ov s} }
(B(\mu, u_\circ^m))^\frac{(2-p)\beta}{q}\!\|u_\circ^m\|_q\,
\|\vp_\circ\|_{\ov s'}\, t^{-\frac{2\beta}{q}}\,, \ea
$$ while the term
$I_\eta$ goes to zero along a subsequence, as $\eta$ tends to zero. Therefore, we get
$$\dy \frac{|(v^{m,\nu}(t),\vf_\circ)|}{  \|\vp_\circ\|_{\ov s'}}
\leq c\, |\O|^{\frac{1}{\ov s} }(B(\mu,
u_\circ^m))^\frac{(2-p)\beta}{q}\!\|u_\circ^m\|_q\,t^{-\frac{2\beta}{q}}
\,\mbox{\,for all }\varphi_\circ\!\in\! C_0^{\infty}(\O),$$ which
implies $$ \|v^{m,\nu}(t)\|_{\ov s}\leq c\,
 |\O|^{\frac{1}{\ov s} }  (B(\mu, u_\circ^m))^\frac{(2-p)\beta}{q}
 \|u_\circ^m\|_q\, t^{-\frac{2\beta}{q}} \,.$$ Since the right-hand side is uniform
with respect to $\ov s$, letting $\ov s\to\infty$, we obtain
$$ \|v^{m,\nu}(t)\|_{\infty}\leq c\,
 (B(\mu, u_\circ^m))^\frac{(2-p)\beta}{q}  \|u_\circ^m\|_q\,t^{-\frac{2\beta}{q}} \,.$$ As in the previous proof, taking into
account suitable strong convergences, we may pass to the limit
firstly as $\nu\to 0$, then as $m\to \infty$, and get the
result.\chiu \vskip0.3cm
 {\bf Proof of Theorem\,\ref{mmtp}} -
 Let $\mu>0$. From Proposition\,\ref{DDAAGG}, $u^\mu\in L^\infty(0,T; V)$ and  $u^\mu_t\in L^2(\ve,T; L^p(\O))$,
  uniformly in $\mu>0$. Moreover, the sequence  $\{u^\mu\}$ converges to the solution
   $u$ of system \eqref{PF}, in suitable norms, as $\mu$ tends to zero (see the proof of
    Theorem\,\ref{existencenomu}).  From the  compact embedding of $V$  in $L^s(\O)$,
    for any $p<s<\frac{np}{n-p}$, and the continuous embedding of   $L^s(\O)$ in $L^p(\O)$,
     using Lemma\,\ref{aubin}, the sequence $\{u^\mu\}$ converges to $u$, as $\mu$ goes
     to zero, strongly in $L^{s}((\ve,T)\times \O)$.
 Hence, $\{u^\mu\}$ converges to $u$ almost everywhere in $t$, $x$, along a subsequence.
  Therefore, along such a subsequence, we find
\be\label{abg0mu}|u(t,x)|\leq |u^\mu(t,x)-u(t,x)|+|u^\mu(t,x)|\leq |u^\mu(t,x)-u(t,x)|
+\|u^{\mu}(t)\|_\infty,\ee
 a.e. in $(t,x)\in (\ve,T)\times \O$. Passing to the limit as $\mu$ goes to zero,
 and using the estimate in Proposition\,\ref{ulinfty}, we easily get
 $$
 \|u(t)\|_\infty\leq \|u_\circ\|_\infty\,.$$
 Therefore we find that  $u$ satisfies \eqref{mmtp1}.
 \par In order to obtain estimate   \eqref{mmtp2}, we can repeat verbatim
  the previous arguments, employing Proposition\,\ref{ul2} in place of
  Proposition\,\ref{ulinfty} to estimate the $L^\infty$-norm of the
  sequence $\{u^\mu\}$ in \eqref{abg0mu}. \chiu

 \section*{\normalsize Appendix}

  \renewcommand{\theequation}{A.\arabic{equation}}
\renewcommand{\thetho}{A.\arabic{tho}}
\renewcommand{\thedefi}{A.\arabic{defi}}
\renewcommand{\therem}{A.\arabic{rem}}
\renewcommand{\theprop}{A.\arabic{prop}}
\renewcommand{\thelemma}{A.\arabic{lemma}}
\renewcommand{\thecoro}{A.\arabic{coro}}
\setcounter{equation}{0} \setcounter{coro}{0} \setcounter{lemma}{0}
\setcounter{defi}{0} \setcounter{prop}{0} \setcounter{rem}{0}
\setcounter{tho}{0}
This appendix is designed for the proof of Proposition\,\ref{existence} and its corollaries.
Firstly we introduce an easy lemma, useful in the sequel.
\begin{lemma}\label{l2lp}{\sl
Let $g(t, x)$ and $F(t,x)$ be such that \be\label{l2lp1}
\|t^{\delta_1}(\mu+|F|^2)^\frac 12\|_{L^\infty(0,T; L^p(\O))}= K_1\,,\quad \|t^{\delta_2}\,(\mu+|F|^2)^\frac{p-2}{4}\, g\|^2_{L^2((0,T)\times\O)}= K_2\,,
\ee
and $t^\delta g(t, x)\in L^2(0,T; L^p(\O))$,
with nonnegative constants $\mu,\delta_1,\delta_2, \delta$ such that
 $\delta:=\frac{2-p}{2}\, \delta_1+\delta_2$. Then, one has
 \be\label{l2lp3}\|t^\delta g(t, x)\|^2_{L^2(0,T; L^p(\O))}\leq K_1^{2-p}\,K_2\,.\ee
}
\end{lemma}
\Pr
Let us write the $L^2(0,T; L^p(\O))$ norm of $t^\delta g(t, x)$ as follows
$$\int_0^t\!\tau^{2\delta} \|g(\tau)\|_p^2\, d\tau=\!
\int_0^t\!\!\tau^{2\delta}\bigg[\int_\O(\mu+|F(\tau,x)|^2)^{\frac{p(p-2)}{4}}
 |g(\tau,x))|^p(\mu+|F(\tau,x)|^2)^{\frac{p(2-p)}{4}}dx\bigg]^\frac 2p\!\!d\tau.$$
By applying H\o lder's inequality with exponents $\frac 2p$ and $\frac{2}{2-p}$, we get
$$\ba{ll}\dy\int_0^t\tau^{2\delta} \|g(\tau)\|_p^2\, d\tau
\\ \dy \leq \int_0^t\!\tau^{\delta_1(2-p)}
\left[\int_\O \tau^{2\delta_2}(\mu+|F(\tau,x)|^2)^\frac{p-2}{2}|g(\tau,x)|^2\,dx\right]\left[
\int_\O(\mu+|F(\tau,x)|^2)^\frac{p}{2} dx\right]^\frac
{2-p}{p} d\tau\ea$$
Hence, from \eqref{l2lp1}  we get
$$\int_0^t\!\!\tau^{2\delta} \|g(\tau)\|_p^2\, d\tau \dy \leq\,\| t^{\delta_1}(\mu+|F|^2)^\frac 12\|_{p}^{2-p}\!\!
 \int_0^t\!\tau^{2\delta_2}\|\,(\mu+|F(\tau)|^2)^\frac {p-2}{4}\,g(\tau)\|_2^2\, d\tau
 \dy \leq
K_1^{2-p}K_2.$$\chiu\vskip0.1cm
 {\bf Proof of Proposition\,\ref{existence}} -  In the sequel we adopt the
  idea introduced by Prodi \cite{prodi} in the context of Navier-Stokes
  equations, where the existence of a solution was proved by the Galerkin
   method with eigenfunctions of the Stokes operator as basis functions.
   Obviously, we replace this basis with the one given by eigenfunctions
   of the Laplace operator. So let $\{a_j\}$ be the eigenfunctions of $\Delta$, and denote
by $\lambda_j$ the corresponding eigenvalues:
$$\ba{ll}\vs1 \dy-\Delta a_j\!\!\!\!&\dy=\lambda_j a_j,\ \mbox{ in }
\O,\\\dy \hskip0.6cm a_j\!\!\!\!&=\dy0,\hskip0.5cm  \mbox{ on }
\po\,. \ea$$  Recall that $a_j\in
W_0^{1,2}(\O)\cap W^{2,2}(\O)$.
 We consider the Galerkin
approximations related to system \eqref{PFepv} of the form
\be\label{G1}v^k(t,x)=\sum_{j=1}^k c_{jk}(t)a_j(x)\ \hskip 0.2cm
k\in \N\,,\ee where the coefficients $c_{jk}$ satisfy the following system
of ordinary differential equations
\be\label{lpsc}\ba{ll}\vspace{0.5ex} \dy \dot
c_{jk}(t)=-\nu\sum_{i=1}^k b_{ji}c_{ik}(t)-
\sum_{i=1}^k d_{ji}c_{ik}(t)=0, \hskip0.2cm  j=1,\ldots,k,\\
\dy c_{jk}(0)=(v_{\circ},a_j),\ea\ee with $b_{ji}=(\nabla a_i,\nabla
a_j)$, $d_{ji}=((\mu+|\nabla (c_{rn}(t)
a_r)|^2)^\frac{(p-2)}{2}\,\nabla a_i,\nabla a_j)$. With this choice
of $c_{jk}(t)$ we impose that $v^{k}(t,x)$ are solutions of the
following system of $k$-differential equations \be\label{ap1}
(v^k_t,a_j)-\nu(\Delta v^k, a_j)+(A_\mu(v^k), a_j)=0\,,
j=1,\cdots,n,\ee with initial conditions $c_{jk}(0)=(v_\circ, a_j)$,
$j=1,\cdots, k$, where $(\cdot,\cdot)$ denotes the standard
$L^2$-inner product. As the right-hand side of \eqref{lpsc} is a
Lipschitz function, due to the assumption $\mu>0$ and
 using Lemma\,\ref{lCarlo2}, the existence of a solution in a time interval
$[0,t_k]$, $t_k\in (0, T]$, follows by standard results on ordinary
differential equations.  The following a priori estimates (see
\eqref{intap2}) will ensure that $t_k=T$, for all $k\in \N$.
\vskip0.2cm {\underline {\it A priori estimates}} - Let us multiply
\eqref{ap1} by $c_{jk}$, by $d c_{jk}/dt$ and sum over $j$. We get,
respectively: \be\label{ap2}
\frac12\frac{d}{dt}\|v^k\|_2^2+\nu\|\nabla v^k\|_2^2+\|(\mu+|\nabla
v^k|^2)^\frac{(p-2)}{4}\nabla v^k\|_2^2=0\,,\ee
 \be\label{ap4a} \|v^k_t\|_2^2+\frac
{\nu}{2}\frac{d}{dt}\|\nabla v^k\|_2^2+\frac
1p\frac{d}{dt}\|(\mu+|\nabla
v^k|^2)^\frac{1}{2}\|_p^p= 0\,.\ee %
Differentiating  \eqref{ap1} with respect to $t$, multiplying by $d
c_{jk}/dt$, then summing over $j$, and observing that
$$\partial_t (a(\mu, v^k)\,\nabla
v^k)= a(\mu,v^k)\,\nabla v^k_t+
(p-2)(\mu+|\nabla
v^k|^2)^\frac{(p-4)}{2}(\nabla
v^k\cdot\nabla v^k_t)\cdot\nabla
v^k\,,$$ we also have
\be\label{ap3aa}
\frac12\frac{d}{dt}\|v^k_t\|_2^2+\nu\|\nabla
v^k_t\|_2^2+(p-1)\|(\mu+|\nabla
v^k|^2)^\frac{(p-2)}{4}\nabla
v^k_t\|_2^2\leq 0\,.\ee The {\it
energy identity} \eqref{ap2}
implies \be\label{intap2}
\|v^k(t)\|_2^2+2\nu\int_0^t\!\|\nabla
v^k(\tau)\|_2^2d\tau+2\int_0^t\!\|a(\mu,
v^k(\tau))^\frac{1}{2}\nabla
v^k(\tau)\|_2^2d\tau=
\|v^k(0)\|_2^2\,,\ee for any $t\in
[0,t_k]$. Since $\|v^k(0)\|_2\leq
\|v_\circ\|_2$, \eqref{intap2}
gives
$$\|v^k(t)\|_2^2=|c_k(t)|^2\leq
\|v_\circ\|_2^2,\ \forall \ t\in
[0,T]\,. $$

Moreover \be\label{nunu}\|\nabla
v^k\|_{L^2(0,T,L^2(\O))}\leq
c(\nu)\,\|v_\circ\|_2\,,\ee and,
using standard arguments\,%
 \footnote{We
are using the inequalities
$$\ba{ll}\vs1\dy\int_{\O}|\nabla u|^p\, dx\!\!&\dy
=\int_{|\nabla u|^2\geq \mu}|\nabla
u|^p\, dx+\int_{|\nabla u|^2 \leq
\mu}|\nabla u |^p\, dx\\& \dy\leq
2^\frac{2-p}{2}\!\int_{\O}
a(\mu,u)\,|\nabla u|^2
dx+\int_{|\nabla u|^2\leq
\mu}\!\!\mu^\frac p2\, dx\leq  c
\int_{\O} a(\mu,u)\, |\nabla
u|^2dx+ \mu^\frac p2\, |\O|.\ea$$}
\be\label{int1ap2}\ba{ll}\dy
\|\nabla
v^k\|^p_{L^p(0,T,L^p(\O))}\leq
c\int_0^T\|a(\mu,
v^k(\tau))^\frac{1}{2}\nabla
v^k(\tau)\|_2^2d\tau +
c(\O,T)\mu^\frac p2 \leq \,
B(\mu,v_\circ)\,,\ea\ee uniformly
in $\nu>0$, with $B$ given by
\eqref{defB}. From this estimate it
also follows that $v^k\in
L^{p}(0,T; V)$ uniformly with
respect to $k$,  and
  $$\|(\mu+|\nabla
v^k|^2)^\frac{(p-2)}{2}\nabla
v^k\|^{p'}_{L^{p'}(0,T;
L^{p'}(\O))}\leq \|\nabla
v^k\|^p_{L^{p}(0,T; L^{p}(\O))}\leq
B(\mu,v_\circ)\,.$$ Hence
\be\label{numA}\|A_\mu(v^k)\|^{p'}_{L^{p'}(0,T;
V')}\leq B(\mu,v_\circ)\,,\ee which
implies $A_\mu(v^k)\in L^{p'}(0,T;
V')$ uniformly with respect to $k$.
Moreover, since $p\in (1,2)$, using
\eqref{nunu} we get
 $$\|(\mu+|\nabla
v^k|^2)^\frac{(p-2)}{2}\nabla
v^k\|_{L^{2}(0,T; L^{2}(\O))}\leq
\|\nabla v^k\|_{L^{2(p-1)}(0,T;
L^{2(p-1)}(\O))}\leq c(\nu)
\,\|v_\circ\|_2^2\,.$$ Hence we
also get
$$\|A_\mu(v^k)\|_{L^{2}(0,T; W^{-1,2}(\O))}\leq
c(\nu) \,\|v_\circ\|_2^2\,.$$ Further $\Delta v^k\in L^{2}(0,T;
W^{-1,2}(\O))$ since, for any $\vp \in  L^{2}(0,T; W_0^{1,2}(\O))$
$$(\Delta v^k,\vp)=-( \nabla v^k,\nabla \vp)\,,$$
and $\nabla v^k\in L^2(0,T; L^2(\O))$ from \eqref{nunu}.
This ensures that $v^k_t\in L^{2}(0,T;
W^{-1,2}(\O))$.
\par
Integrating \eqref{ap4a} from $0$ to $t$, we get
\be\label{a13a}\ba{ll}\dy \int_0^t \|v^k_\tau(\tau)\|_2^2d\tau+\frac
{\nu}{2}\,\|\nabla v^k(t)\|_2^2+\frac 1p\,\|(\mu+|\nabla
v^k(t)|^2)^\frac{1}{2}\|_p^p\\\hskip1cm  \dy\leq  \frac
{\nu}{2}\|\nabla v^k(0)\|_2^2+ \frac 1p\|(\mu+|\nabla
v^k(0)|^2)^\frac{1}{2}\|_p^p
\dy\leq  c\nu\, \|\nabla v_\circ\|_2^2+ c(\O,T)\,\mu^\frac p2, \ea\ee %
where, observing that $c_{lk}(0)=(v_\circ,
a_l)=\frac{1}{\lambda_l}(\nabla v_\circ,\nabla a_l)=\frac{(\nabla
v_\circ,\nabla a_l)}{\|\nabla a_l\|_2^2}$, we have used
$$\|\nabla v^k(0)\|_2\leq \|\nabla v_\circ\|_2\,. $$%
 \par
Integrating inequality (\ref{ap3aa}) from $0$ to $t$, we get
\be\label{a16aa}\ba{ll}\dy\!\!\|v^k_t(t)\|_2^2\!+2\nu\!\!\int_0^t\!\!\!\!\|\nabla
v^k_\tau(\tau)\|_2^2d\tau+2(p-1)\!\int_0^t\!\!\!\|a(\mu,v^k(\tau))^\frac{1}{2}\nabla
v^k_\tau(\tau)\|_2^2d\tau\!\leq \!\|v^k_t(0)\|_2^2. \ea\ee Let us
estimate the right-hand side.  Multiplying \eqref{ap1} by $d
c_{jk}/dt$ and summing over $j$ we have $$\|v^k_t\|_2^2\leq
\nu\,\|\Delta v^k\|_2\|v^n_t\|_2+  \|A_\mu(v^k)\|_2\|v^k_t\|_2\,\leq
(\nu+c(\mu))\,\|\Delta v^k\|_2\|v^k_t\|_2\,,$$ since
$$\ba{ll}\dy \vs1
\|A_\mu(v^k)\|_2^2&\!\dy= \left\|\frac{\Delta v^k}{(\mu+|\nabla
v^k|^2)}{\null_{\frac {2-p}{2}}} +(p-2) \frac{\nabla
v^k\cdot\nabla\nabla v^k \cdot\nabla v^k}{\ \ (\mu+|\nabla
v^k|^2)^\frac{(4-p)}{2}} \right\|_2^2\\& \dy \leq
2\,\mu^{p-2}\,\|\Delta v^k\|_2+2(2-p)\,\mu^{p-2}\,\|D^2
v^k\|_2^2\leq c(\mu)\, \|\Delta v^k\|_2^2\,.\ea $$ Observing that,
as $c_{lk}(0)=(v_\circ, a_l)=\frac{(\Delta v_\circ,\Delta
a_l)}{\|\Delta a_l\|_2^2}$, then $\|\Delta v^k(0)\|_2\leq \|\Delta
v_\circ\|_2$, we get \be\label{vnt0}\|v^k_t(0)\|_2^2\leq
(\nu+c(\mu)) \|\Delta v_\circ\|_2\,.\ee Finally, by using estimates \eqref{a13a}, \eqref{a16aa}
and \eqref{vnt0}, we can apply Lemma\,\ref{l2lp}, with $g=\nabla
v^k_t$,  $F=\nabla v^k$, $\delta_1=\delta_2=\delta=0$, and obtain
$$\int_0^t\,\|\nabla
v^k_\tau(\tau)\|_p^2\, d\tau
 \dy \leq  (\nu+c(\mu)) \|v_\circ\|_{2,2}^2
 \left( c\nu\|\nabla v_\circ\|_2^2+ c(\O,T)\mu^\frac
 p2\right)^\frac{2-p}{p}.$$
 \vskip0.2cm {\underline {\it  Passage to the limit}} -
Using the above estimates we can extract a subsequence, still
denoted by $\{v^k\}$, such that, in the limit as $k$ tends to
$\infty$, uniformly in $\nu>0$, $$ v^k\, \rightharpoonup \, v
\textrm{ in } L^{\infty}(0,T; L^2(\O)) \ \textrm{weakly}-*\,;$$
\be\label{con2} v^k\, \rightharpoonup \, v \textrm{ in }
L^{p}(0,T;V)\ \textrm{weakly}\,; \ee \be\label{con3} \,v^k_t\,
\rightharpoonup \, \,v_t \textrm{ in } L^{\infty}(0,T; L^2(\O))\
\textrm{weakly}-*\,; \ee $$ \,v^k(t)\, \rightharpoonup \, \,\xi\
\textrm{ in }  L^2(\O)\ \textrm{weakly}\,; $$ $$\nabla v^k_t\,
\rightharpoonup \,  \nabla v_t \textrm{ in } L^{2}(0,T; L^p(\O))\
\textrm{weakly}\,; $$ $$A_\mu(v^k)\rightharpoonup \chi\ \textrm{in }
L^{p'}(0,T; V')\ \textrm{weakly}.$$
 The convergence \eqref{con3} implies convergence in
$L^{p'}(0,T;L^2(\O))$, hence \be\label{con3a}\,v^k_t\,
\rightharpoonup \, \,v_t \textrm{ in } L^{p'}(0,T; V')\
\textrm{weakly},\ee
and
$\xi=v(t)$. Further, from \eqref{con2} and \eqref{con3a}, $v\in C(0,T; L^2(\O))$.
 Moreover, non-uniformly in
$\nu$, we also have $$ v^k\, \rightharpoonup \, v \textrm{ in }
L^{2}(0,T; W^{1,2}_0(\O))\ \textrm{weakly}\,,$$
$$A_\mu(v^k)\rightharpoonup \chi\ \textrm{in } L^{2}(0,T;
W^{-1,2}(\O))\ \textrm{weakly},$$ $$v^k_t\rightharpoonup v_t\
\textrm{in } L^{2}(0,T; W^{-1,2}(\O))\ \textrm{weakly},$$
 Following usual arguments from the theory
of monotone operators (see \cite{lions}, Chapter 2, Sec. 1), one
shows that the non-linear part $A_\mu(v^k)$ actually (weakly)
converges to $A_\mu(v)$, as $k$ tends to $\infty$ and that the limit $v$ is a weak solution of
system \eqref{PFepv}. The regularities stated for $v$ follows from
the analogous regularities of $\{v^k\}$ and the lower semi-continuity of
the norm for the weak convergence. Finally, the uniqueness follows
from the monotonicity of $A_\mu$.\chiu

\vskip0.3cm
  {\bf Proof of Corollary\,\ref{existencecor1}} - We argue as
  in the proof of Proposition\,\ref{existence},
  without exploiting the $C_0^\infty$-  regularity of the initial data.
  So, let us consider \eqref{ap4a}.  Multiplying by $t$, a simple
computation gives
 \be\label{a12}\ba{ll}\dy\vs1
t\,\|v^k_t\|_2^2+\frac {\nu}{2}\frac{d}{dt} (t\,\|\nabla v^k\|_2^2)+\frac
1p\frac{d}{dt}(t\,\|(\mu+|\nabla v^k|^2)^\frac{1}{2}\|_p^p)\\
\dy \hfill= \frac
{\nu}{2}\|\nabla v^k\|_2^2+ \frac 1p\|(\mu+|\nabla
v^k|^2)^\frac{1}{2}\|_p^p .\ea\ee Integrating inequality (\ref{a12})
from $0$ to $t$ and estimating the right-hand side of \eqref{a12}
via inequalities \eqref{intap2} and \eqref{int1ap2}, we get
\be\label{a13}\ba{ll}\dy \vs1\int_0^t \tau
\|v^k_\tau(\tau)\|_2^2d\tau+\frac {\nu}{2}\,t\,\|\nabla
v^k(t)\|_2^2+\frac tp\,\|(\mu+|\nabla v^k(t)|^2)^\frac{1}{2}\|_p^p\\
\dy\hskip1cm = \frac {\nu}{2}\int_0^t\|\nabla v^k(\tau)\|_2^2d\tau+ \frac
1p\int_0^t\|(\mu+|\nabla v^k(\tau)|^2)^\frac{1}{2}\|_p^pd\tau
\dy\leq  c\,B(\mu, v_\circ), \ea\ee %
with $B$ given in \eqref{defB}.
Hence, in particular, we obtain
 \be\label{os4} \|t^\frac 12\, v^k_t\|_{L^2(0,T;L^2(\O))}+ t\,\|\nabla
v^k(t)\|_p^p\leq c\,B(\mu,
v_\circ)\,.\ee \par Let us consider
\eqref{ap3aa}. Multiplication by
$t^2$ gives \be\label{a15}
\frac12\frac{d}{dt}(t^2\|v^k_t\|_2^2)+\nu
\,t^2\|\nabla
v^k_t\|_2^2+(p-1)t^2\|(\mu+|\nabla
v^k|^2)^\frac{(p-2)}{4}\nabla
v^k_t\|_2^2\leq
t\,\|\,v^k_t\|_2^2\,.\ee
Integrating inequality (\ref{a15})
from $0$ to $t$, and then
estimating the right-hand side
 via  \eqref{a13}, we get
\be\label{a16}\ba{ll}\dy
t^2\|\,v^k_t(t)\|_2^2+2\nu\int_0^t\!\!\tau^2\|\nabla
v^k_\tau(\tau)\|_2^2\,d\tau\\ \dy
+2(p-1)\int_0^t\!\!\tau^2\|(\mu+|\nabla
v^k(\tau)|^2)^\frac{(p-2)}{4}\nabla
v^k_\tau(\tau)\|_2^2\,d\tau\leq
2\int_0^t\!\!\tau\,\|\,v^k_\tau(\tau)\|_2^2d\tau\leq
c\, B,\ea\ee which ensures that
$t\,\nabla v^k_t\in
L^{2}(0,T;L^2(\O))$, non-uniformly
in $\nu>0$. Finally, by using
estimates  \eqref{a13} and
\eqref{a16}, we can apply
Lemma\,\ref{l2lp}, with $g=\nabla
v^k_t$,  $F=\nabla v^k$,
$\delta_1=\frac 1p$, $\delta_2=1$,
hence $\delta=\frac{p+2}{2p}$, and
obtain
\be\label{a18}\ba{ll}\dy\int_0^t\tau^{\frac
{p+2}{p}}\,\|\nabla
v^k_\tau(\tau)\|_p^2\, d\tau \leq
c\,B^\frac 2p(\mu, v_\circ).\ea\ee
 The previous bounds \eqref{os4}, \eqref{a16} and \eqref{a18}
ensure that, up to subsequences, in the limit of $k\to \infty$: $$t^\frac 1p \nabla v^k
\rightharpoonup t^\frac 1p \nabla v\ \textrm{ in }\ L^\infty(0,T;
L^p(\O))\ \textrm{weakly}-*,$$
$$t v^k_t
\rightharpoonup t v_t\  \textrm{ in }\ L^\infty(0,T; L^2(\O)) \ \textrm{weakly}-*, $$
$$tA_\mu(v^k_t)\rightharpoonup t\chi_t\ \textrm{in }
L^{2}(0,T; W^{-1,2}(\O))\ \textrm{weakly},$$
$$t^\frac{p+2}{2p}\nabla v^k_t
\rightharpoonup t^\frac{p+2}{2p}\nabla v_t\  \textrm{ in } L^2(0,T;
L^p(\O))\ \textrm{weakly},$$ and completes the proof.\chiu  \vskip0.3cm
   {\bf Proof of Corollary\,\ref{pregcor}}   -  Let us consider the Galerkin approximations
\eqref{G1}. Let us multiply \eqref{ap1} by $\l_j c_{jk}$ and sum
over $j$. By observing that
$(v^k_t,\sum_{j=1}^k\l_jc_{jk}a_j)=(v^k_t,-\Delta v^k)=
\frac{d}{dt}\|\nabla v^k\|_2^2$, recalling Lemma\,\ref{LL1} and using
Cauchy's inequality, we get \be\label{ap5}\ba{ll}\vs1\dy
 \frac12\frac{d}{dt}\|\nabla
v^k\|_2^2+\nu\|\Delta v^k\|_2^2+ \,\|(\mu+|\nabla
v^k|^2)^\frac{(p-2)}{4}\Delta v^k\|_2^2\\ \dy \leq
\left((2-p)\,C_1+\frac{\delta}{2}\right) \,\|(\mu+|\nabla
v^k|^2)^\frac{(p-2)}{4}\Delta v^k\|_2^2+ \frac{C}{2\delta}\,\left(
\|\nabla v^k\|_p^p+ \mu^\frac p2|\O|\right)
 \,,\ea\ee
 for any $\delta>0$. Hence, since $p>p_\circ\geq \frac 32$,
 choosing $\delta=1-(2-p)C_1\equiv\overline C(p)$ we get the inequality
\be\label{appri6}\frac{d}{dt}\|\nabla
v^k\|_2^2+\overline C(p)
\,\|(\mu+|\nabla
v^k|^2)^\frac{(p-2)}{4}\Delta
v^k\|_2^2\leq C \|\nabla
v^k\|_p^p+\, C\,\mu^\frac p2
|\O|\,. \ee By using that
$v_\circ\in W_0^{1,2}(\O)$,
integrating \eqref{appri6} from $0$
to $t$, we have
$$\ba{ll}\dy \vs1\|\nabla
v^k(t)\|_2^2+\overline C(p)
\int_0^t\,\|(\mu+|\nabla
v^k(\tau)|^2)^\frac{(p-2)}{4}\Delta
v^k(\tau)\|_2^2d\tau \\
\hskip3cm\dy \leq \|\nabla
v_\circ\|_2^2+ C \int_0^t(\|\nabla
v^k(\tau)\|_p^p+\, C\,\mu^\frac p2
|\O|)d\tau\,,\ea $$ which ensures
that  $\nabla v^k\in
L^{\infty}(0,T; L^2(\O))$ uniformly
with respect to  $k$, and, applying
ones again Lemma\,\ref{LL1},
  $(\mu+|\nabla v^k|^2)^\frac{(p-2)}{4}D^2 v^k\in L^2(\O)$ uniformly
with respect to  $k$.
  Moreover, by using estimates \eqref{a13a}, \eqref{a16aa}
and \eqref{vnt0}, and  employing Lemma\,\ref{l2lp},
  with $g=D^2v^k$,  $F=\nabla v^k$, $\delta_1=\delta_2=\delta=0$, we obtain $$\ba{ll}\dy\vs1 \int_0^t\|D^2
v^k(\tau)\|_p^2\,d\tau
 \\ \dy\hfill \leq \,\|(\mu+|\nabla
v^k|^2)^\frac{1}{2}\|_{L^\infty(0,T; L^p(\O))}^{2-p}\, \int_0^t \|\,a(\mu,
v^k(\tau))^\frac{1}{2}\,D^2 v^k(\tau)\|_2^2\, d\tau
 \dy \leq
C(B(\mu, v_\circ),T).\ea$$ If we do not exploit the $W_0^{1,2}$-
regularity of the initial data, but just its $L^2$-integrability, we
can argue as follows. Throughout this proof, we denote $B(\mu,
v_\circ)$ just by $B$. Multiplication of \eqref{appri6} by
$t^\alpha$, for some $\alpha\geq 1$ that will be specified later,
gives
$$\ba{ll}\vs1\dy \frac{d}{dt}(t^\alpha\|\nabla v^k\|_2^2)+\overline
C(p)t^\alpha\,\|(\mu+|\nabla v^k|^2)^\frac{(p-2)}{4}\Delta
v^k\|_2^2\\\hfill\dy \leq \alpha\, t^{\alpha-1}\,\|\nabla v^k\|_2^2+
C t^\alpha \|\nabla v^k\|_p^p+\, C t^\alpha \mu^\frac p2|\O|\,.\ea$$
An integration from $0$ to $t$ gives \be\label{ap7}\ba{ll}\vs1\dy
t^\alpha\|\nabla v^k(t)\|_2^2+\,\overline C(p)\int_0^t
\tau^\alpha\|(\mu+|\nabla v^k(\tau)|^2)^\frac{(p-2)}{4}\Delta
v^k(\tau)\|_2^2\,d\tau\\\hskip1cm\dy \leq
\alpha\int_0^t\tau^{\alpha-1} \|\nabla v^k(\tau)\|_2^2\,d\tau +
C\int_0^t \left(\tau^\alpha \|\nabla
v^k(\tau)\|_p^p+\,\tau^\alpha\mu^\frac p2|\O|\right)\,d\tau\,.\ea\ee
For the first integral on the right-hand side we argue as follows.
Observing that $p>\frac{2n}{n+2}$, by using a Gagliardo-Nirenberg
inequality we get
$$ \|\nabla v^k\|_2^2\leq  C\|D^2 v^k\|_p^{2a}\|\nabla v^k\|_p^{2(1-a)},
 \quad a=n\mbox{\large$\left(\frac 1p-\frac 12\right)$}\in (0,1)\,.$$ Observing that
$$\|D^2 v^k\|_p^{2}\leq \|(\mu+|\nabla v^k|^2)^\frac{(p-2)}{4}D^2 v^k\|_2^2\, \|(\mu+|\nabla
v^k|^2)^\frac{1}{2}\|_p^{2-p}\,,$$ we have
$$\ba{ll}\vs1\dy \|\nabla v^k(\tau)\|_2^2\leq C \, \|(\mu+|\nabla v^k|^2)^\frac{(p-2)}{4}
D^2 v^k\|_2^{2a}\,\|(\mu+|\nabla
v^k|^2)^\frac{1}{2}\|_p^{(2-p)a}\,\|\nabla v^k\|_p^{2(1-a)}\,,
 \ea$$
and, by Cauchy's inequality, for any $\delta>0$, we get
 \be\label{ago1}\ba{ll}\dy\vs1  \tau^{\alpha-1}\|\nabla v^k(\tau)\|_2^2\\\vs1\dy\leq
 \frac{\delta}{2}
 \tau^{\alpha}\|a(\mu, v^k)^\frac{1}{2} D^2 v^k\|_2^{2} +
\frac{C}{2\delta}\tau^{\alpha-\frac{1}{1-a}} \|(\mu+|\nabla
v^k|^2)^\frac{1}{2}\|_p^{\frac{(2-p)a}{1-a}}\|\nabla
v^k\|_p^{2}\\\dy \leq \frac{\delta}{2} \tau^{\alpha}\|a(\mu,
v^k)^\frac{1}{2} D^2
v^k\|_2^{2}\!+\frac{C}{2\delta}\,\tau^{\alpha-\frac{1}{1-a}}
\mu^\frac{(2-p)a}{2(1-a)}\|\nabla
v^k\|_p^2+\frac{C}{\delta}\,\tau^{\alpha-\frac{1}{1-a}} \|\nabla
v^k\|_p^\frac{2-pa}{1-a}.\ea\ee The first term on the right-hand
side can be estimated using Lemma\,\ref{LL1}. Let us integrate the
last two terms on the right-hand side of \eqref{ago1} from $0$ and
$t$. Since
$$\int_0^t\,\tau^{\alpha-\frac{1}{1-a}} \|\nabla
v^k(\tau)\|_p^\frac{2-pa}{1-a}\,d\tau=\int_0^t\,\tau^{\alpha-\frac{1}{1-a}}
\|\nabla v^k(\tau)\|_p^{\frac{2-p}{1-a}}\|\nabla v^k(\tau)\|_p^p\,d\tau\,, $$ by
choosing $\alpha$ in such a way that
$$\alpha=\mbox{\large$\frac{1}{1-a}+\frac{2-p}{p(1-a)}$},
$$
hence $\alpha$ as in \eqref{alpha}, and then
using \eqref{os4} and \eqref{int1ap2}, we get
$$\int_0^t\,\tau^{\alpha-\frac{1}{1-a}} \|\nabla
v^k(\tau)\|_p^\frac{2-pa}{1-a}\,d\tau\leq
B^\frac{2-p}{p(1-a)}\int_0^t\|\nabla v^k(\tau)\|_p^p\,d\tau\leq
B^\frac{2-ap}{p(1-a)}.$$ Fixed $\alpha$, in a similar
way we also get
$$\ba{ll}\vs1\dy\int_0^t\tau^{\alpha-\frac{1}{1-a}}
\|\nabla v^k(\tau)\|_p^2\, d\tau=\int_0^t\tau^{\frac{2-p}{p(1-a)}} |\nabla
v^k(\tau)\|_p^{2-p}\,
 \|\nabla v^k(\tau)\|_p^{p}\, d\tau\\ \dy \hskip2cm\leq
B^{\frac{2-p}{p}}t^{\frac{2-p}{p}\frac{a}{1-a}}\int_0^t\|\nabla
v^k(\tau)\|_p^{p}\, d\tau\leq
B^{\frac{2}{p}}t^{\frac{2-p}{p}\frac{a}{1-a}}.\ea
$$
Therefore, integrating \eqref{ago1} from $0$ to $t$ and using the
previous estimates for the terms on the right-hand side, we obtain
$$\ba{ll}\dy\vs1\int_0^t\tau^{\alpha-1} \|\nabla
v^k(\tau)\|_2^2\,d\tau
 \leq C_1\frac{\delta}{2} \int_0^t\,\tau^{\alpha}\|a(\mu,v^k(\tau))^\frac{1}{2}
\Delta v^k(\tau)\|_2^{2}\,d\tau
\\\hskip1cm \dy +C\frac{\delta}{2}\int_0^t\!\!\tau^\alpha(\|\nabla v(\tau)\|_p^p+
\mu^\frac p2 |\O|)\, d\tau
 + \frac{C}{2\delta}\,
\mu^\frac{(2-p)a}{2(1-a)}\,B^\frac 2p\,
t^{\frac{2-p}{p}\frac{a}{1-a}}+\frac{C}{2\delta}\,B^\frac{2-ap}{p(1-a)}\,.\ea$$
Inserting this estimate in \eqref{ap7}, then choosing
$\delta=\frac{\overline C(p)}{\alpha C_1}$ and using \eqref{a13}, we
arrive at $$\ba{ll}\vs1\dy t^\alpha \|\nabla
v^k(t)\|_2^2+\,\frac{\overline C(p)}{2}\int_0^t\tau^\alpha \|a(\mu,
v^k)^\frac{1}{2}\,\Delta v^k(\tau)\|_2^2\,d\tau\\\hfill\dy \leq C\,
\mu^\frac{(2-p)a}{2(1-a)}\,B^\frac 2p\,
t^{\frac{2-p}{p}\frac{a}{1-a}}+\frac{C}{2\delta}\,B^\frac{2-ap}{p(1-a)}
+c B t^\alpha+c Bt^{\alpha+1}\,.\ea$$ Observing that
$\frac{2-p}{p}\frac{a}{1-a}<\alpha$, and defining $\beta_1(p)=\frac
12(\alpha-\frac{2-p}{p}\frac{a}{1-a})=\frac{1}{2}(\alpha-
\frac{n(2-p)^2}{p(2p-n(2-p))})$, as in \eqref{beta12} the above estimate shows in
particular that $t^{\beta_1(p)}\, \nabla v^k \in L^\infty(0,T;
L^2(\O))$. This estimate, together with $\nabla v^k\in C(\ve,T;
L^p(\O))$, which follows from Corollary\,\ref{existencecor1},  gives
$\ t^{\beta_1(p)}\, \nabla  v^k \in C_w(0,T; L^2(\O)).$ The strong
continuity can be obtained as follows. From \eqref{appri6} we have,
for any $t>s>0$,
 \be\label{fwc}\|\nabla v^k(t)\|_2^2\leq C \int_s^t
(\,\|\nabla v^k(\tau)\|_p^p+\, C\,\mu^\frac p2 |\O|)\, d\tau+ \|\nabla v^k(s)\|_2^2\,.\ee
By using the identity
$$\|\nabla v^k(t)-\nabla v^k(s)\|_2^2 = \|\nabla v^k(t)\|_2^2+
\|\nabla v^k(s)\|_2^2-2(\nabla v^k(t),\nabla v^k(s))\,,$$
then estimate \eqref{fwc} and the  the weak continuity of
$\nabla v^k(t)$ in $L^2(\O)$ we get the result.
\par Finally, by  applying Lemma\,\ref{l2lp}, with $g=D^2
v^k$,  $F=\nabla v^k$, $\delta_1=\frac 1p$, $\delta_2=\frac{\alpha}{2}$, hence $\delta=\frac{p+2}{2p}$, and obtain
 $$\ba{ll}\dy\int_0^t\tau^{\alpha+\frac{2-p}{p}} \|D^2
v^k(\tau)\|_p^2\, d\tau
\\ \dy \leq t^\frac{2-p}{p}\,\|(\mu+|\nabla
v^k|^2)^\frac{1}{2}\|_{p}^{2-p}\, \int_0^t\tau^\alpha \|\,a(\mu,
v^k(\tau))^\frac{1}{2}\,D^2 v^k(\tau)\|_2^2\, d\tau
 \dy \leq
C(B,T).\ea$$ The previous bounds,
all uniform with respect to $k$,
ensure the weak-* convergence of a
 subsequence of $\{t^{\beta_1(p)}\, \nabla v^k\}$ in the space
$L^\infty(0,T; L^2(\O))$ and, recalling the expression of $\beta_2(p)$ given in \eqref{beta12}, the weak convergence of a
 subsequence of $\{t^{\beta_2(p)}\, D^2v^k\}$ in the space $L^2(0,T; L^p(\O))$, as $k\to\infty$,
 uniformly in $\nu$, $\mu$. From Proposition\,\ref{existence},
 we get that the limit solution $v$ of \eqref{PFepv}
satisfies  $t^{\beta_1(p)}\,\nabla v\in L^\infty(0,T; L^2(\O))$ and
$t^{\beta_2(p)}\,\, v\in L^2(0,T; W^{2,p}(\O))$.
 \chiu
\vskip0.1cm
 {\bf Acknowledgment} -
The paper is
 performed under the
auspices of GNFM-INdAM.


\begin{thebibliography}{20}

\bibitem{AcMin} E. Acerbi and G. Mingione, {\it Gradient estimates for a class of parabolic systems}, Duke Math. J., {\bf  136} (2007), 285--320.

\bibitem{AMS} E. Acerbi, G. Mingione and G.A. Seregin,
{\it Regularity results for parabolic systems related to a class of non-Newtonian fluids},
Ann. Inst. H. Poincar\'{e} Anal. Non Lin\'{e}aire, {\bf  21} (2004), 25--60.


\bibitem{BDV} H. Beir\~ao da Veiga, {\it Singular parabolic $p$-Laplacian systems under non-smooth external forces.  Regularity up to the boundary}, online arXiv:1206.1808v1 [math.AP].

\bibitem{CDB} Y.Z. Chen and E. DiBenedetto,  {\it  Boundary estimates for solutions of nonlinear
degenerate parabolic systems}, J. Reine Angew. Math. 395 (1989),
102--131.

\bibitem{choe} H. Choe, {\it H\"{o}lder continuity of solutions of certain degenerate
parabolic systems}, Non-linear Anal., {\bf  8} (1992), 235--243.

\bibitem{CG}
F. Crispo and C. R. Grisanti, {\it On the existence,
 uniqueness and $C^{1,\gamma}(\overline\Omega)\cap W^{2,2}(\Omega)$
  regularity for a class of shear-thinning fluids},  J. Math. Fluid Mech., {\bf 10}
   (2008), 455--487.
   \bibitem{CMellittico} F. Crispo and P.
   Maremonti, {\it Higher regularity
of solutions to the singular
\mbox{$p\,$}-Laplacean system},
submitted.


\bibitem{DB} E. DiBenedetto, {\it Degenerate Parabolic Equations}, Universitext,
New York, Springer- Verlag. XV, 1993.

\bibitem{DBF} E. DiBenedetto and A. Friedman, {\it H\"{o}lder estimates for non-linear
degenerate parabolic systems}, J. Reine Angew. Math., {\bf 357}
(1985), 1--22.

\bibitem{DBH} E. DiBenedetto and M.A. Herrero, {\it
Nonnegative solutions of the evolution p-Laplacian equation. Initial traces and Cauchy problem when $1<p<2$},
Arch. Rational Mech. Anal., {\bf  111} (1990), 225--290.
\bibitem{DBKV} E. DiBenedetto, Y.C. Kwong and V.
Vespri, {\it Local space analicity of solutions of certain singular
parabolic equations}, Indiana Univ. Math. J., {\bf 40} (1991),
741--765.

\bibitem{DBUV} E. DiBenedetto, J.M. Urbano and V.
Vespri, {\it Current Issues on Singular and Degenerate Evolution
Equations}, in Handbook of Differential Equations - Evolutionary
Equations, vol. I, Elsevier North Holland 2004, 169--286.

\bibitem{DMSt} F. Duzaar, G.
Mingione and K. Steffen, {\it Parabolic Systems with Polynomial
Growth and Regularity}, Mem. Amer. Math. Soc., {\bf  214} (2011).

\bibitem{KL} J. Kinnunen and J.L. Lewis, {\it Higher integrability for parabolic systems of $p$-Laplacian type},
 Duke Math. J., {\bf  102} (2000), 253--271.

\bibitem{L} O.A. Ladyzhenskaya,
{\it The mathematical theory of viscous incompressible flow}, Gordon
and Breach, 1968.

\bibitem{LU} O.A. Ladyzhenskaya and N.N. Ural'tseva, {\it Linear and quasilinear elliptic equations},
 Academic Press, New York-London 1968.

  \bibitem{LSU} O.A. Ladyzhenskaya, V.A. Solonnikov, N.N. Ural'ceva,
{\it Linear and Quasi-linear equations of Parabolic Type},
Providence American Mathematical Society, 1968.

 \bibitem{LVP} C. Leone,  A. Verde and G. Pisante,  {\it
Higher integrability results for non smooth parabolic systems: the
subquadratic case}, Discrete Contin. Dynam. Syst. B, {\bf 11}
(2009), 177--190.


\bibitem{lions} J.-L. Lions,  {\it Quelques m\'{e}thodes de r\'{e}solution
des probl\`{e}mes aux limites non lin\'{e}aires}, Dunod;
Gauthier-Villars, Paris, 1969.

\bibitem{MS} P. Maremonti and V.A.
Solonnikov, {\it On nonstationary
Stokes problem in exterior
domains}, Ann. Scuola Norm. Sup. Pisa Cl. Sci., { \bf 24} (1997), 395--449.

\bibitem{prodi} G. Prodi, {\it Teoremi di tipo locale per il sistema
di Navier-Stokes e stabilit\`{a} delle soluzioni stazionarie}, Rend.
Sem. Mat. Padova, {\bf 32} (1962), 374-397.

\bibitem{Sch1} C. Scheven, {\it  Nonlinear Calder\'{o}n-Zygmund theory for parabolic systems with subquadratic
growth}, J. Evol.  Equ., {\bf 10} (2010), 597--622.

\bibitem{Sch2} C. Scheven, {\it Regularity for subquadratic parabolic systems: higher integrability and dimension
estimates}, Proc. Roy. Soc. Edinburgh, {\bf  140 A} (2010), 1269--1308.

\bibitem{show} R.E. Showalter, {\it Monotone Operators in Banach Space and Nonlinear Partial Differential Equations}, American Mathematical Society 1997.

\bibitem{Sol77}  V. A. Solonnikov, {\it Estimates for solutions of nonstationary Navier-Stokes equations},
 J. Soviet Math., {\bf 8} (1977), 467--528.


\end{thebibliography}
\end{document}